\documentclass[a4paper,reqno]{amsart}
\usepackage{a4,wasysym}
\usepackage{amsmath,graphicx, amssymb,color,ulem,bbm,geometry}
\usepackage{epsf, subfigure, verbatim}
\usepackage{float}
\usepackage{epsfig}
\usepackage{enumitem}

\newtheorem{theorem}{Theorem}[section]
\newtheorem{corollary}[theorem]{Corollary}
\newtheorem{lemma}[theorem]{Lemma}

\newtheorem{definition}[theorem]{Definition}

\newcommand{\ud}{\mathrm{d}}

\newcommand{\X}{\mathcal{X}}

\newcommand{\be}{\begin{equation}}
\newcommand{\ee}{\end{equation}}

\newcommand{\eofproof}{\hfill$\Box$ } 

\makeatletter
\newcommand{\mylabel}[2]{#2\def\@currentlabel{#2}\label{#1}}
\makeatother

\begin{document}

\newgeometry{left=1.75cm, right=1.75cm, top=2.3cm, bottom=2cm}

\numberwithin{equation}{section}

\title[On the formulation of size-structured consumer resource models ]{On the formulation of size-structured consumer resource models 
(with special attention for the principle of linearised stability)}

\author[C. Barril]{Carles Barril}
\address{Carles Barril, Department of Mathematics, Universitat Aut\`{o}noma de Barcelona, Bellaterra, 08193, Spain}
\email{carlesbarril@mat.uab.cat}

\author[\`{A}. Calsina]{\`{A}ngel Calsina}
\address{\`{A}ngel Calsina, Department of Mathematics, Universitat Aut\`{o}noma de Barcelona and Centre de Recerca Matem\`{a}tica,  Bellaterra, 08193, Spain}
\email{acalsina@mat.uab.cat}

\author[O. Diekmann]{Odo Diekmann}
\address{Odo Diekmann, Department of Mathematics, University of Utrecht, Budapestlaan 6, PO Box 80010, 3508
TA, Utrecht, The Netherlands}
\email{O.Diekmann@uu.nl}

\author[J. Z. Farkas]{J\'{o}zsef Z. Farkas}
\address{J\'{o}zsef Z. Farkas, Division of Computing Science and Mathematics, University of Stirling, Stirling, FK9 4LA, United Kingdom }
\email{jozsef.farkas@stir.ac.uk}

\date{\today}

\begin{abstract}

To describe the dynamics of a size-structured population and its unstructured resource, we formulate bookkeeping equations in two different ways. The first, called the PDE formulation, is rather standard. It employs a first order partial differential equation, with a non-local boundary condition, for the size-density of the consumer, coupled to an ordinary differential equation for the resource concentration. The second is called the DELAY formulation and employs a renewal equation for the population level birth rate of the consumer, coupled to a delay differential equation for the (history of the) resource concentration. With each of the two formulations we associate a constructively defined semigroup of nonlinear solution operators.

The two semigroups are intertwined by a non-invertible operator. In this paper we delineate in what sense the two semigroups are equivalent. In particular, we i) identify conditions on both the model ingredients and the choice of state space that guarantee that the intertwining operator is surjective, ii) focus on large time behaviour and iii) consider full orbits, i.e., orbits defined for time running from $-\infty$ to $+\infty$.

Conceptually, the PDE formulation is by far the most natural one. It has, however, the technical drawback that the solution operators are not differentiable, precluding rigorous linearisation. (The underlying reason for the lack of differentiability is exactly the same as in the case of state-dependent delay equations: we need to differentiate with respect to a quantity that appears as argument of a function that may not be differentiable.) For the delay formulation, one can (under certain conditions concerning the model ingredients) prove the differentiability of the solution operators and establish the Principle of Linearised Stability. Next the ‘equivalence’ of the two formulations yields a rather indirect proof of this principle for the PDE formulation.

\end{abstract}
\maketitle

\section{Introduction}

When formulating a structured population model, one starts by specifying the so-called i-states, i.e. the states that individuals can have. Next one specifies the relevant aspects of the external world and captures them by variables that describe the environmental condition as experienced by the individuals. The model specification concerns the behaviour of an individual, in particular its i-state development, survival and reproduction (not only the number of offspring, but also their state-at-birth, has to be specified), as well as the impact on (i.e. feedback to) the environmental condition. Often one specifies how rates depend on the current i-state and the prevailing environmental condition, but as demonstrated in \cite{DGM1,DGM2}, an attractive and more general (meaning that i-state development is not necessarily deterministic) alternative is to introduce a composite model ingredient for i-state development and survival over a non-infinitesimal period of time, in combination with a description of cumulative reproductive output over such a period. 

Once the i-level model is complete, it is a bookkeeping exercise  to lift it to the p-level (population level). At this stage a choice has to be made whether to work with measures or restrict to densities. Interpretation is a good guide when working with measures, while one needs greater care when working with densities. Yet, when there is no need to work with measures, modellers usually work with densities. As explained in detail in \cite{MD}, the temporal change of a density is described by a first order PDE with, as a general rule, non-local terms, essentially since the i-state of offspring is far from the i-state of the mother, in general. Integration along characteristics yields quasi-explicit expressions for the solution when an initial density is prescribed and non-local terms are replaced by `given' functions. A consistency requirement then leads to an equation that has to be solved in order to obtain a true solution. Often the equation originating from the consistency requirement is a renewal equation, and existence and uniqueness of solutions is readily established. Inserting the solution in the formula for the density at later times, one obtains a dynamical system describing how the density changes over time. 

More recently it has been argued (see e.g. \cite{DGMNR}) that one can take the renewal equation itself as the starting point for the definition of a dynamical system. The renewal equation is interpreted as a rule for extending a function of time towards the future on the basis of the (assumed to be) known past, so as a delay equation. By shifting along the extended 
function, i.e.,  by updating the history, as for instance described in \eqref{DE-4} below, one obtains a dynamical system \cite{D2}. 

A natural question arises: how do the PDE and the delay dynamical systems relate to each other? One would not like to obtain different dynamics for one and the same i-level model, when population bookkeeping is done in a different manner. In \cite{CDF} we investigated this question for a linear model in which the distribution of i-state-at-birth is described by a density, a so-called distributed states at birth model. Here we investigate the question for a nonlinear consumer-resource model, 
in which all of the consumer newborns have the same i-state (a single state at birth model).  That is, we focus here on a class of structured consumer-resource models, which describe the interaction and population dynamics of a size-structured consumer and its unstructured resource. On the one hand there is significant intrinsic mathematical interest in these nonlinear models, as they pose analytical and computational challenges (see for example \cite{Breda,DGMNR,FH,Grabosch1989,Grabosch1990,Grabosch1991,SBDGV,HT}). On the other hand, for particular choices of the model ingredients, they are also used to investigate or demonstrate the richness of the dynamical behaviour of, for instance, a size-structured population of {\em Daphnia} feeding on algae (see e.g. \cite{Adoteye,DeRoos, deRP}).

The specific model we consider has as i-state a one-dimensional quantity that describes the size of an individual. It is usually denoted by $x$ or $\xi$, and it takes values from the interval $[x_b,\infty)$, where $x_b$ is the size at birth. 
We assume that individuals cannot shrink and their growth is deterministic, which in particular means that two individuals of the same size experiencing the same environment will have the same size throughout their life. The environmental condition (that consumers experience) is determined 
by the resource (food) concentration and it is denoted by $S$. The model ingredients are as follows:
\begin{itemize}
\item the growth rate $g(x,S)$, which we assume to be positive,
\item the death rate $\mu(x,S)$, assumed to be non-negative,
\item the reproduction rate $\beta(x,S)$, assumed to be non-negative,
\item the resource consumption rate $\gamma(x,S)$, assumed to be non-negative,
\item the rate $f(S)$ of change of $S$ in the absence of consumers.
\end{itemize}
We will impose regularity conditions on the model ingredients later on. Dynamic energy budget theory, see \cite{K, wikiDEBT},  provides relations between, on the one hand, the per capita gain through ingestion $\gamma$ and, on the other hand, the per capita expenditure for metabolism, growth $g$ and reproduction $\beta$. Such relations are quite important, but as they are irrelevant for the analysis of this paper, we shall not dwell on them here.

This paper has two main aims. First, we want to describe precisely how the two dynamical systems corresponding to the same population model relate to each other. In some sense, the set of initial conditions for the delay equation formulation is `bigger' than the corresponding set for the PDE formulation. But the difference is inessential in that it does not affect the population birth rate b and the resource concentration S, the two variables that are constructed for $t>0$ from the initial condition and the model ingredients. Once b and S are constructively defined, there is an explicit expression for each of the two semigroups. So the delay formalism has a certain redundancy. By concentrating on the essential information we establish the asymptotic equivalence of the two nonlinear semigroups.

Secondly, we want to prove the Principle of Linearised Stability for steady states. A direct verification for the PDE formulation is impossible, simply since the semigroup operators are NOT differentiable (this observation which, as far as we know, has not been made before, explains why the literature so far does not contain results about linearised stability for size structured models with variable growth rate). The reason is that the initial population density is not only reduced by death, but also translated, by growth, over a variable distance. And if the initial density is not absolutely continuous, there is no differentiable dependence on this distance.

In case of delay equations, bookkeeping is based on ‘time since’, i.e., age, and translation has a fixed rather than a variable speed, so the difficulty disappears. The most straightforward path to linearised stability is by way of a linearisation of the delay equations (but see \cite{DK}). For infinite delay and differentiable equations, \cite{DG} provides a proof of the principle. As far as we know, a rigorous proof of differentiability of the equations corresponding to this kind of model has not been given before (note that \cite{DGMNR} and \cite{DGMNR-r} show that even formal linearisation is not that easy when the growth rate is allowed to have a discontinuity). Here we determine conditions on the model ingredients that allow us to prove the differentiability of the equations. (For sure these conditions are too restrictive and hopefully future work will relax them.)
     
Once the principle is established in the delay equation setting, the asymptotic equivalence results allow us to transfer it to the PDE setting. Thus we circumvented the problem of non-differentiability. In fact our results justify stability conclusions in the PDE setting based on information about the roots of a characteristic equation obtained by formal linearisation of the PDE (the ‘formal’ being that one differentiates an unbounded operator without any attention for its domain). 

We are going to impose certain conditions on the model ingredients to prove the equivalence between the two formulations and the differentiability of the nonlinear operators appearing in the delay formulation. Here we introduce shorthand notations for the various hypotheses. Below, the letter $h$ may refer to any of the model ingredients, i.e. $h\in\{f,g,\mu,\beta,\gamma\}$.
\begin{itemize}
\item[\mylabel{Hh-1}{H1$_h$}] The function $h$ is globally Lipschitz continuous.
\item[\mylabel{Hh-2}{H2$_h$}] The function $h$ is bounded from above.
\item[\mylabel{Hh-3}{H3$_h$}] The function $h$ is bounded from below by a positive constant.
\item[\mylabel{Hh-4}{H4$_h$}] The function $h$ is continuously differentiable.
\item[\mylabel{Hh-5}{H5$_h$}] The function $h$ is differentiable and $Dh$ is globally Lipschitz continuous, i.e. for all $y_0$ and $y$ in the domain of $h$
\[
\|Dh(y_0)-Dh(y)\|\leq L \|y_0-y\|.
\]
where $L$ is a constant independent of $y_0$ and $y$.
\item[\mylabel{Hginf}{H$_{g_\infty}$}]
    There exists $g_\infty>0$ and $\bar{x}>x_b$, such that the function $g$ satisfies $g(x,S)=g_\infty$ for all $x\geq\bar{x}$, and $S\geq 0$.
    
     The final hypothesis involves the functions $\mu$ and $g$.
\item[\mylabel{Hs}{H$_{s}$}]  $\mu(x,S)=\hat{\mu}+\tilde\mu(x,S)$ with $\hat{\mu}>0$ and $\tilde\mu$ (not necessarily positive), such that
\begin{equation*}
    |\tilde\mu(x,S)|\leq \sigma(x)g(x,S),
\end{equation*}
    for a positive function $\sigma$ that is integrable over $[x_b,\infty)$, i.e.,
\begin{equation*}
    \int_{x_b}^{\infty}\sigma(x) dx <\infty .
\end{equation*}
\end{itemize}

The structure of the paper is as follows. In Section 2 we elaborate the formulation of the model in terms of a first order PDE for the size-density of the consumer. By integration along characteristics we derive a renewal equation for the population level birth rate. After a brief discussion of existence and uniqueness, we define the PDE dynamical system.
 In Section 3 we explain the DELAY formulation and next introduce the Banach space of weighted histories that serves as the state space for the corresponding dynamical system.
 In Section 4 we specify a map $\mathcal{L}$ that maps the DELAY state space to the PDE state space. We group the elements of the DELAY state space that are mapped to the same element of the PDE state space into an equivalence class and we show that, under certain conditions on the model ingredients, $\mathcal{L}$ has a pseudo-inverse. For full orbits, i.e., orbits that go back in time to $-\infty$, we establish a one-to-one relationship.
 Section 5 is devoted to steady states. We show that these are characterised by one equation in one unknown and that stability in the PDE setting is equivalent to stability in the DELAY setting. In the DELAY setting we derive a characteristic equation and next formulate the Principle of Linearised Stability.
 The final Section 6 is devoted to some concluding remarks. Technical proofs are provided in three appendices.

\section{The PDE formulation}\label{sectPDEformulation}
Let $\kappa_0\geq0$ and $L^1_{\kappa_0}:=L^1_{\kappa_0}([x_b,\infty);\mathbb{R})$ the space of integrable functions with the weighted norm 
\begin{equation*}
\|n\|_{L^1_{\kappa_0}}=\int_{x_b}^\infty |n(x)|e^{\kappa_0 x}\,\ud x,
\end{equation*}
so that $L^1_0$ is the space of integrable functions, whereas for $\kappa_0>0$ the space $L^1_{\kappa_0}$ consists of a proper subset of the space of integrable functions. Indeed, notice that if $n\in L^1_{\kappa_0}$ then the number in the tail beyond $x$, i.e. $\int_x^\infty n$, decays exponentially at a rate $-\kappa_0$ as $x$ tends to $\infty$.  We are interested both in $\kappa_0=0$ and $\kappa_0>0$, for technical reasons that will be explained in subsequent sections. We denote the positive cone by $L^1_{\kappa_0,+}$.

Let $n(t,\cdot)\geq 0$ denote the density of the size distribution of the consumer population and let $S(t)$ denote the resource concentration  at time $t\ge 0$. Our aim is to determine $n(t,\cdot)\in L^1_{\kappa_0}$ and $S(t)\ge 0$, for $t\ge 0$ from the initial conditions
\begin{equation}\label{PDE-initial}
\begin{aligned}
n(0,x)&=n_0(x)\geq 0,\quad & n_0\in L^1_{\kappa_0}, \\ 
S(0)&=S_0\geq 0,\quad & 
\end{aligned}
\end{equation}
by solving the system of equations
\begin{equation}\label{PDE}
\begin{aligned}
\frac{\partial n}{\partial t}(t,x)+\frac{\partial}{\partial x}\left(g(x,S(t))n(t,x)\right)&=-\mu(x,S(t))n(t,x), \\ 
g\left(x_b,S(t)\right)n(t,x_b)&=\int_{x_b}^\infty\beta(\xi,S(t))n(t,\xi)\,\ud\xi, \\
\frac{\ud S}{\ud t}(t)&=f(S(t))-\int_{x_b}^\infty\gamma(\xi,S(t))n(t,\xi)\,\ud \xi.
\end{aligned}
\end{equation}

The initial conditions are points in the product space $L^1_{\kappa_0,+}\times\mathbb{R}_+$, which we refer to as the space of population densities and environmental conditions. As the norm of a point $(n_0,S_0)$ in this space we choose
\begin{equation}\label{norm-densities}
\|(n_0,S_0)\|_{\kappa_0}=\|n_0\|_{L^1_{\kappa_0}}+|S_0|.
\end{equation}

The PDE in (\ref{PDE}) describes changes in the density $n$ due to growth and death 
of individuals, while the boundary condition (second equation in \eqref{PDE}) determines changes in $n$ due to reproduction. The function $f$ determines the intrinsic dynamics of the resource, that is, $f$ determines changes of the resource population which are not due to the consumer population. Given resource concentration $S$, a consumer of size $x$ consumes on 
average per unit of time $\gamma(x,S)$ units of resource, and so the second term at the right hand-side of the last equation in \eqref{PDE} captures the change in $S$ due to consumption.

In a first step towards a constructive definition of the solution of \eqref{PDE-initial}-\eqref{PDE}, we pretend that both $S(t)$ and the p-level birth rate 
\begin{equation}\label{birth-rate}
b(t)=\int_{x_b}^\infty\beta(x,S(t))n(t,x)\,\ud x
\end{equation}
are given functions of time, for time in appropriate intervals. With this 
in mind we introduce
\begin{equation}\label{X-s-def}
\begin{aligned}
X_S(t,s,\xi)= &\, \text{size of an individual at time}\,\, t\,\, \text{given the individual has size}\,\, \xi\,\, \text{at time}\,\, s \\ 
&\, \text{and given}\,\, S\,\, \text{between times}\,\, s\,\, \text{and}\,\, t,
\end{aligned}
\end{equation}
while noting that often, but not always, we have $t>s$; and
 
\begin{equation}\label{F-s-def}
\begin{aligned}
\mathcal{F}_S(t,s,\xi)= &\, \text{probability that an individual of size}\,\,\xi\,\, \text{at time}\,\, s\,\, \text{is still alive at time}\,\, t>s, \\ 
&\, \text{given}\,\, S\,\, \text{between times}\,\, s\,\, \text{and}\,\, t.
\end{aligned}
\end{equation}
More formally we define
\begin{equation}\label{X-s-def-2}
X_S(t,s,\xi):=x(t),
\end{equation}
where $x$ is the unique solution of 
\begin{equation*}
\dot{x}(\tau)=g(x(\tau),S(\tau)),\quad x(s)=\xi,
\end{equation*}
and
\begin{equation}\label{F-s-def-2}
\mathcal{F}_S(t,s,\xi):=\exp\left(-\int_s^t\mu(x(\tau),S(\tau))\,\ud\tau\right).
\end{equation}


In addition we introduce
\begin{equation}\label{tau-def}
\begin{aligned}
T_S(x,\xi,s) =& \,\text{time at which size equals}\,\,x,\,\, \text{given size equals}\,\, \xi\,\, \text{at time}\,\,s\,\, \text{and}\\
&\,\text{given}\,\,\tau\mapsto S(\tau),\,\, \text{ with}\,\,\tau\,\,\text{between}\,\,s\,\,\text{and}\,\,T_S,
\end{aligned}
\end{equation}
$T_S$ and $X_S$ are inverse functions in the sense that
\begin{equation}\label{TXinverse-eq}
T_S (X_S (t,s,\xi),\xi,s) = t
\end{equation}
and
\begin{equation}\label{XTinverse-eq}
X_S(t,T_S(\xi,x,t),\xi)=x.
\end{equation}

So in particular, 

\begin{equation}\label{tau-eq}
X_S(t,T_S(x_b,x,t),x_b)=x.
\end{equation}

Using these relations it follows that $t(x):=T_S(x,\xi,s)$ is the unique solution of
\begin{equation*}
t'(x)=\frac{1}{g(x,S(t(x)))},\quad t(\xi)=s.
\end{equation*}

\begin{center}
\begin{figure}[h!]
\begin{centering}
\includegraphics[width=5.5in]{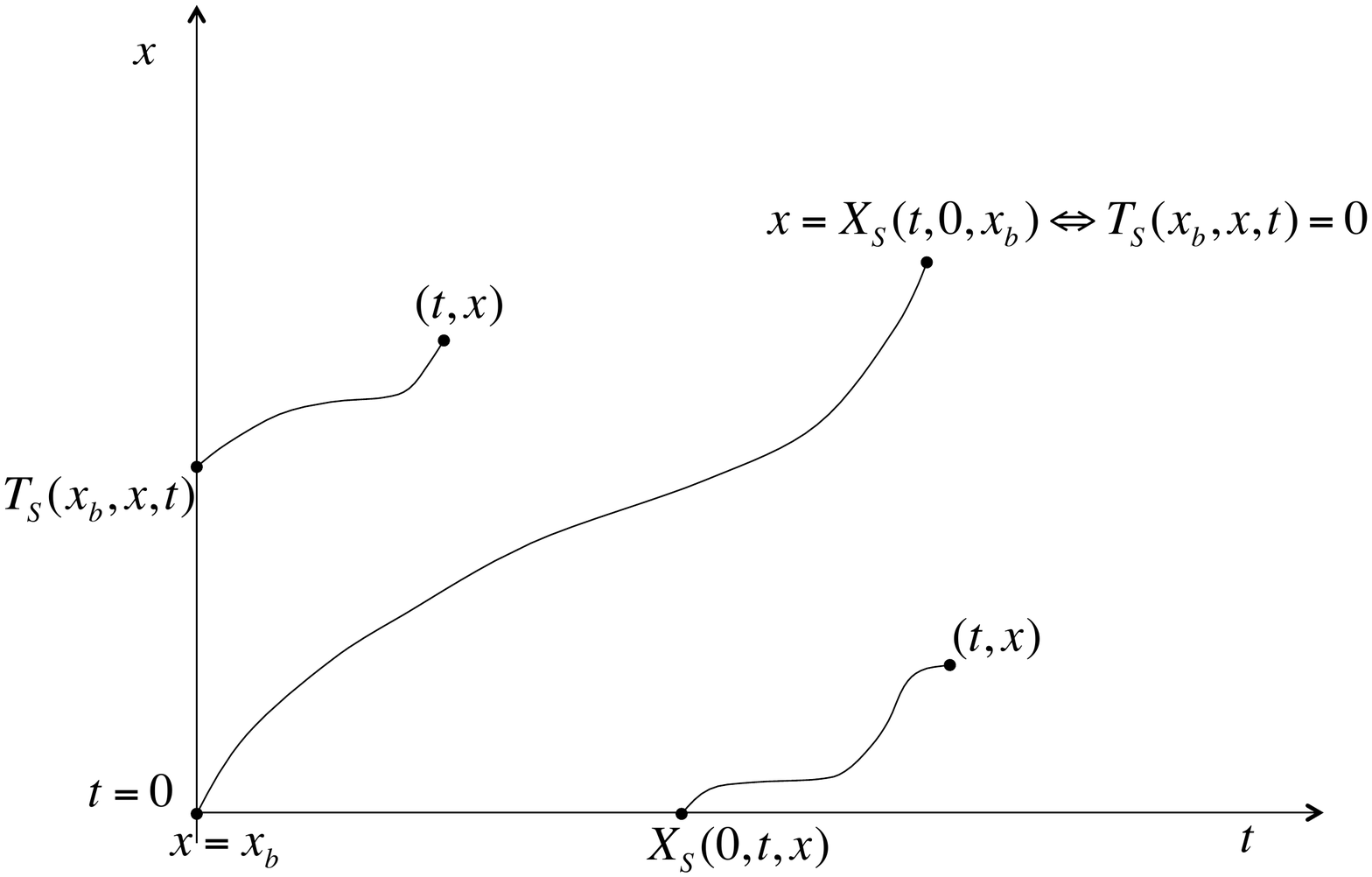}
\end{centering}
\caption{\label{fig1} For two different time-size combinations $(t,x)$, the solution of the equation $\dot{x}(\tau)=g(x(\tau),S(\tau))$ is followed backward in time. In one case the trajectory hits the time-axis at $T_S(x_b,x,t)$, in the other case it hits the size-axis at $X_S(0,t,x)$. The boundary curve between the two cases is the curve $x=X_S(t,0,x_b)$ parametrised by $t$. This curve can be also characterised by the equation 
$T_S(x_b,x,t)=0$. }
\end{figure}
\end{center}

The key aspects of integration along characteristics are formulated in the next Lemma and illustrated in Figure \ref{fig1}. 

\begin{lemma}
For given $S$ and $b$, the solution of 
\begin{equation}\label{PDE-2}
\begin{aligned}
\frac{\partial n}{\partial t}(t,x)+\frac{\partial}{\partial x}\left(g(x,S(t))n(t,x)\right)&=-\mu(x,S(t))n(t,x) \\ 
g\left(x_b,S(t)\right)n(t,x_b)&=b(t) \\
n(0,x)&= n_0(x)
\end{aligned}
\end{equation}
is given explicitly by
\begin{equation}\label{char-sol}
n(t,x)=\left\{\begin{aligned} & n_0(X_S(0,t,x))\mathcal{F}_S(t,0,X_S(0,t,x))D_3X_S(0,t,x),& \text{when}\quad X_S(0,t,x)>x_b 
\\ & b(T_S(x_b,x,t))\mathcal{F}_S(t,T_S(x_b,x,t),x_b)(-D_2\,T_S(x_b,x,t)), &  \text{when}\quad\, T_S(x_b,x,t)>0 \end{aligned} \right.
\end{equation}
while $n(t,x)$ is not specified when $T_S(x_b,x,t)=0$, or equivalently $X_S(0,t,x)=x_b$.
\end{lemma}

Deliberately, we have not specified beforehand what we do mean by a solution of \eqref{PDE-2} and therefore we cannot provide a proof of this lemma. What we provide instead is a more fundamental conservation principle that leads directly to \eqref{char-sol}. In this view, \eqref{PDE-2} is a concise symbolic infinitesimal representation of the conservation principle. 

Let $N(t,\cdot)$ be the cumulative size distribution (i.e., the representation of the measure with density $n(t,\cdot)$ by an NBV function, normalized to be zero at $x=x_b$). The ‘explicit’ formula
\begin{equation}\label{NBVrep-eq1}
N(t,x) = \int_{T_S(x_b,x,t)}^t b(\tau) \mathcal{F}_S(t,\tau,x_b)\,\ud\tau   \qquad\text{for}\quad x < X_S(t,0,x_b),    
\end{equation}
expresses that $N(t,x)$ consists, for small $x$, of individuals born after time zero who have survived till time $t$ and have not (yet) grown beyond size $x$. For larger $x$ values we have
\begin{equation}\label{NBVrep-eq2}
N(t,x) = \int_0^t  b(\tau) \mathcal{F}_S(t,\tau,x_b)\,\ud\tau+\int_{x_b}^{X_S(0,t,x)}\mathcal{F}_S(t,0,\xi) N_0\,(\ud\xi)\qquad\text{for}\quad x > X_S(t,0,x_b),
\end{equation}
expressing that all individuals born after time zero and before time t are included, provided they survive, and, in addition, those individuals already present at time zero who survived and have not (yet) grown beyond $x$. So here $N_0$ denotes the initial cumulative distribution and $b$ the birth rate (considered to be given). These formulas follow directly from the interpretation and cover more general initial conditions (i.e., we do not need to restrict to $N_0$ being absolutely continuous). See \cite{DGM2,DGM1} for an analysis of general structured population models along these lines.  

Formula (\ref{char-sol}) is obtained by formal differentiation of (\ref{NBVrep-eq1}) and (\ref{NBVrep-eq2}) with respect to $x$. 

It is straightforward to derive, by integration of (\ref{PDE-2}) with respect to $x$, an equation for $N$ and to solve this equation by integration along characteristics as taught in PDE text books. This leads to (\ref{NBVrep-eq1})-(\ref{NBVrep-eq2}). Alternatively one can define measure solutions via duality (note that (\ref{PDE-2}) is the Kolmogorov forward equation for the density and that one can ‘lift’ it to measures via the corresponding backward equation and duality), see \cite{G}. In \cite{SBDGV} the cumulative formulation is used as the starting point for numerical work.

Substituting \eqref{char-sol} into the right-hand side of \eqref{birth-rate} we find an equation for $b(t)$ (in terms of the trajectory of $S$ between 0 and $t$). Using instead of $x$ the variable $\xi$ corresponding to either the birth size $x_b$ or the size at $t=0$ and featuring in our derivation of \eqref{char-sol} above, we obtain the renewal equation
\begin{equation}\label{RE}
b_S(t)=\int_0^t\beta_S(t,t-a)b_S(t-a)\,\ud a+h_S(t),
\end{equation}
where
\begin{equation}\label{RE-2}
\beta_S(t,s):=\tilde{\beta}_S(t,s,x_b)
\end{equation}
with
\begin{equation}\label{RE-3}
\tilde{\beta}_S(t,s,\xi):=\beta\left(X_S(t,s,\xi),S(t)\right)\mathcal{F}_S(t,s,\xi)
\end{equation}
being the expected contribution to the population birth rate at time $t$ of an individual that has size $\xi$ at time $s<t$, and where
\begin{equation}\label{f-def}
h_S(t):= \int_{x_b}^\infty n_0(\xi)\tilde{\beta}_S(t,0,\xi)\,\ud\xi,
\end{equation}
sums the expected contributions to the population birth rate at time $t$ of individuals that existed at time zero and had size $\xi$ at that time.

We can solve \eqref{RE} by {\it generation expansion}, i.e. by successive 
approximation. If we inductively define
\begin{equation}\label{beta-1}
\beta_S^{k+1}(t,s):=\int_s^t\beta_S^k(t,\tau)\beta_S(\tau,s)\,\ud\tau,\quad k\ge 1,
\end{equation}
then $\beta_S^k(t,s)$ is the rate at which the k-th generation offspring of an individual born at time $s$ is produced at time $t>s$. The so-called clan kernel $\beta^c_S$ is obtained by summing over all generations:
\begin{equation}\label{beta-2}
\beta_S^c(t,s):=\sum_{k=1}^\infty\beta_S^k(t,s).
\end{equation}
See \cite{DGM1} for more details on the resolvent representation of the solution of the renewal equation above. We have:
\begin{lemma}
The unique solution of the renewal equation \eqref{RE} is given by
\begin{equation}\label{RE-sol}
b_S(t)=h_S(t)+\int_0^t\beta_S^c(t,\tau)h_S(\tau)\,\ud\tau.
\end{equation}
\end{lemma}

It remains to determine the resource concentration as a function of time, 
i.e. it remains to solve the initial value problem 
\begin{equation}\label{initial-resource}
\begin{aligned}
\frac{\ud S}{\ud t}(t)= & f(S(t))-\int_0^t\gamma_S(t,t-a)b_S(t-a)\,\ud a-k_S(t), \\ 
S(0)= & S_0,
\end{aligned}
\end{equation}
where
\begin{equation}\label{gamma-1}
\gamma_S(t,s):=\tilde{\gamma}_S(t,s,x_b),
\end{equation}
with
\begin{equation}\label{gamma-2}
\tilde{\gamma}_S(t,s,\xi):=\gamma\left(X_S(t,s,\xi),S(t)\right)\mathcal{F}_S(t,s,\xi)
\end{equation}
being the expected rate at which an individual which has size $\xi$ at time $s$ consumes substrate at time $t$, and where
\begin{equation}\label{ks-def}
k_S(t):=\int_{x_b}^\infty n_0(\xi)\tilde{\gamma}_S(t,0,\xi)\,\ud\xi
\end{equation}
is the substrate consumption rate at time $t$ due to the individuals that were already present at time zero and survived till time $t$.
\begin{lemma} \label{lemma2.4}
Assume that H1$_h$ for $h\in\{f,g,\mu,\beta,\gamma\}$, H2$_\beta$ and H2$_\gamma$ hold. Then \eqref{initial-resource} has a unique global solution (for positive times). 
\end{lemma}

{\bf Sketch of the proof.}
One can show that for all $\hat{t}>0$ the function $\mathcal{V}:C([0,\hat{t}],\mathbb{R})\rightarrow C([0,\hat{t}],\mathbb{R})$ defined as
\begin{equation}\label{volterraeq}
\mathcal{V}(S)(t)= S_0 + \int_0^t \left( f(S(\tau))-\int_0^\tau \gamma_S(\tau,\tau-a)b_S(\tau-a)\,\ud a - k_S(\tau)\right)d\tau,
\end{equation}
is a contraction when using a suitable norm on $C([0,\hat{t}],\mathbb{R})$, which is equivalent to the supremum norm. For full  details of the proof see Appendix A.
\eofproof

\begin{theorem}\label{Theorem 2.4}
Assume that H1$_h$ for $h\in\{f,g,\mu,\beta,\gamma\}$, H2$_\beta$ and H2$_\gamma$ hold. Then problem \eqref{PDE-initial}-\eqref{PDE} has a unique global solution. Moreover, the family of continuous (solution) operators 
\begin{equation*}
\mathcal{T}_{PDE}(t)\,:\,L^1_{\kappa_0}\times\mathbb{R}_+\to L^1_{\kappa_0}\times\mathbb{R}_+,
\end{equation*}
defined by 
\begin{equation}\label{sg-pde}
\mathcal{T}_{PDE}(t)\begin{pmatrix} n_0 \\ S_0 \end{pmatrix}=\begin{pmatrix} n(t,\cdot) \\ S(t) \end{pmatrix},
\end{equation}
is strongly continuous and has the semigroup property.
\end{theorem}

As \cite{CS} contains a detailed proof of the corresponding result for a slightly different model, we refrain from providing the proof of this theorem. We refer to \cite{HT} for interesting considerations concerning the relation between model assumptions and well-posedness.

It is not too difficult to characterise the steady states of \eqref{PDE-initial}-\eqref{PDE}, i.e. solutions $(n_*,S_*)$ satisfying 
\begin{equation}\label{ss=pde}
\mathcal{T}_{PDE}(t)\begin{pmatrix} n_* \\ S_* \end{pmatrix}=\begin{pmatrix} n_* \\ S_* \end{pmatrix},\quad \forall\,t\ge 0,
\end{equation}
see also \cite{DGM}. In Section 5 we are going to discuss in more detail the existence of positive steady states, and we also refer the interested 
reader to \cite{CF2,FH}, where the steady state problem for a (more general) consumer-resource model was addressed. 

Once $S$ is constructed, the solution operator for the size density is given by (\ref{char-sol}). If we fix $n_0$, but vary $S(0)$, there are changes in the argument $X_S(0,t,x)$ of $n_0$. So differentiability with respect to $S(0)$ requires that $n_0$ is absolutely continuous, which, in general, it is not. We conclude that the nonlinear solution operators are not differentiable and that, consequently, we cannot linearise at a steady state. This is reminiscent of the lack of differentiability for state-dependent delay equations, \cite{HKWW}. And the underlying reason is identical: we have to differentiate a quantity that occurs as argument of a function that need not be differentiable. We conclude that within the PDE framework we cannot even formulate the Principle of Linearised Stability, let alone prove it.

\section{The delay equation formulation}

If the resource concentration $S$ is (considered to be) known for negative times, then in Figure 1 there is no need to stop when a trajectory hits 
the size-axis, one can instead continue the trajectory backwards in time until it hits the size$=x_b$ axis. The earlier interpretation of $\beta_S$ and $\gamma_S$ then directly leads to the system of equations
\begin{equation}\label{DE-1}
\begin{aligned}
b(t)=& \int_0^\infty \beta_S(t,t-a)b(t-a)\,\ud a, \\
\frac{\ud S}{\ud t}(t)=& f(S(t))-\int_0^\infty\gamma_S(t,t-a)b(t-a)\,\ud a,
\end{aligned}
\end{equation}
to which we add an `initial' condition in the form of a prescribed history for both $b$ and $S$ at a certain time, say zero (system \eqref{DE-1} is time-translation invariant in a sense that can be made precise, so when 
choosing zero as the time at which we prescribe the history we do not lose generality):
\begin{equation}\label{DE-initial}
\begin{aligned}
b(\theta)=& \phi(\theta) \\
S(\theta)=& \psi(\theta)
\end{aligned},\quad \theta\le 0,
\end{equation}
for non-negative functions $\phi$ and $\psi$ that we consider as given. The function $\phi$ should be locally integrable, while $\psi$ should be continuous, and shortly we will add conditions concerning the growth of $\phi$ and $\psi$ as $\theta\to -\infty$. 

The same argument that led to the second expression for $n(t,x)$ in \eqref{char-sol} now yields that, given \eqref{DE-initial}, the size-distribution at time zero is given by
\begin{equation}\label{DE-2}
n_0(x)=\phi(T_\psi(x_b,x,0))\mathcal{F}_\psi(0,T_\psi(x_b,x,0),x_b)\left(-D_2\,T_\psi(x_b,x,0)\right),
\end{equation}
and clearly the resource concentration at time zero is given by
\begin{equation}\label{DE-3}
S_0=\psi(0).
\end{equation}
It follows that the issue of existence and uniqueness of solutions is already covered by our discussion in Section 2: $b(t)$ for $t>0$ is defined by \eqref{RE-sol} with $h_S$ defined by \eqref{f-def} with $n_0$ defined by \eqref{DE-2}, and $S(t)$ is for $t>0$ the solution of \eqref{initial-resource} with $S_0$ defined by \eqref{DE-3}.

We view \eqref{DE-1} as a rule for extending functions of time towards the future on the basis of the known past. A dynamical system is obtained by translation along the extended pair of functions, i.e. by updating the history
\begin{equation}\label{DE-4}
\mathcal{T}_{DE}(t)\begin{pmatrix} \phi \\ \psi \end{pmatrix}:=\begin{pmatrix} b_t \\ S_t \end{pmatrix}
\end{equation}
where we employ the usual notational convention
\begin{equation}\label{DE-5}
q_t(\theta):=q(t+\theta),\quad \theta\le 0.
\end{equation}

As we want to make use of the results of \cite{DG}, we want $\phi$ and $\psi$ to belong to weighted function spaces. In order to have a natural choice for the weight in the $\phi$ component, we impose the following assumption:

\noindent There exist $\hat\mu>0$, $0<c\leq 1$ and $C\ge 1$, such that
\begin{equation}\label{weight-ass1}
c\,e^{-\hat\mu(t-s)}\le \mathcal{F}_S(t,s,x_b)\le C\,e^{-\hat\mu(t-s)},\quad \forall\, t\ge s.
\end{equation}
Here we allow for $C>1$ in order to incorporate models where small individuals do not suffer any mortality, so that the survival probability equals one until some time after birth. 

Let us show that assumption \eqref{weight-ass1} holds under hypothesis $\text{H}_{s}$.

\begin{lemma} Assume $\text{H}_{s}$. Then \eqref{weight-ass1} holds.
\end{lemma}

{\bf Proof.} Applying $\text{H}_{s}$ in  \eqref{F-s-def-2} we have: 
\begin{equation*}
\mathcal{F}_S(t,s,x_b)=e^{-\hat\mu (t-s)}\exp\left\{-\int_s^t\tilde\mu(x(\tau),S(\tau))\,\ud\tau\right\}=e^{-\hat\mu (t-s)}\exp\left\{-\int_{x_b}^{x(t)}\tilde{\mu}(\xi,S(x^{-1}(\xi)))\,\frac{1}{x'(x^{-1}(\xi))}\,\ud\xi\right\},
\end{equation*}
where $x(\tau)=X_S(\tau,s,x_b)$ and $x^{-1}(\xi)$, the inverse of $x(\tau)$, can be interpreted as the time at which an individual, with size $x_b$ at time $s$, reaches size $\xi$ (notice that since $x(\tau)$ depends on the function $S$, so does $x^{-1}(\xi)$). Now, since $x'(x^{-1}(\xi))=g(\xi,S(x^{-1}(\xi)))$, from H$_{s}$ it follows that
\begin{equation*}
-\int_{x_b}^{\infty}\sigma(\xi)\,\ud\xi\leq -\int_{x_b}^{x(t)}\frac{\tilde{\mu}(\xi,S(x^{-1}(\xi)))}{g(\xi,S(x^{-1}(\xi)))}\,\ud\xi \leq \int_{x_b}^{\infty}\sigma(\xi)\,\ud\xi,
\end{equation*}
so that, if we define $c$ and $C$ by
\begin{equation*}
\log(c):=-\int_{x_b}^{\infty}\sigma(\xi)\,\ud\xi,\quad \log(C):=\int_{x_b}^{\infty}\sigma(\xi)\,\ud\xi,
\end{equation*}
we have
\begin{equation*}
 c \leq \exp\left\{-\int_{x_b}^{x(t)}\tilde{\mu}(\xi,S(x^{-1}(\xi)))\,\frac{1}{x'(x^{-1}(\xi))}\ud\xi\right\} \leq  C,
\end{equation*}
and then \eqref{weight-ass1} follows.
\eofproof

The growth condition on $\phi$ is now expressed by the requirement that
\begin{equation}\label{weight-ass2}
\int_0^\infty\left|\phi(-a)\right|e^{-\mu_0 a}\,\ud a<\infty,
\end{equation}
where $\mu_0\in(0,\hat\mu]$, so that the total number of individuals is bounded whatever the birth history is. The limit case $\mu_0=\hat\mu$ corresponds to the biggest set of birth histories, but (as explained in Appendix \ref{appendix-diffF}) that choice could lead to delay equations that fail to be differentiable. This is why an exponent $\mu_0$ smaller or equal than $\hat\mu$ is considered. A key point in all of this is that constant functions should be in the state space, since 
we  want to consider steady states and their stability.
To summarise: we want $\mu_0$ to be 
\begin{itemize}
    \item[-] positive, in order to:
    \begin{enumerate}
        \item[i.] include steady states in the state space,
        \item[ii.] guarantee that, as explained in \cite{DG}, the essential spectrum is restricted to the open left half of the complex plane;
    \end{enumerate}
    \item[-] less than $\hat{\mu}$ in order to avoid that a large birth rate in the distant past can lead to a population of infinite size. 
\end{itemize} 
This still leaves some freedom, which we are going to exploit in Appendix \ref{appendix-diffF} when discussing the differentiability of the nonlinear delay equations.

Concerning $\psi$, there does not seem to be a natural growth condition, so we rather arbitrarily work with the exponent $\mu_0$ for $\psi$ too, by 
assuming that
\begin{equation}\label{weight-ass3}
\lim_{a\to\infty}\left|\psi(-a)\right|e^{-\mu_0 a}=0.
\end{equation}

\begin{definition}\label{def-phase-space-DE}
Let $\mu_0\in (0,\hat{\mu}]$ with $\hat{\mu}$ being the mortality rate for old individuals as in (\ref{weight-ass1}). Then let
\begin{align*}
||\phi||_{1,\mu_0}:=& \int_0^\infty|\phi(-a)|e^{-\mu_0a}\,\ud a, \\
||\psi||_{\infty,\mu_0}:=& \sup_{0\le a<\infty}\left\{|\psi(-a)|e^{-\mu_0 a}\right\}, \\
\left|\left|\begin{pmatrix}
\phi \\ \psi \end{pmatrix} \right|\right|_{\mu_0}:=& \,||\phi||_{1,\mu_0}+||\psi||_{\infty,\mu_0}.
\end{align*}
 The (positive cone of) the state space is defined as
\begin{equation*}
\mathcal{X}_{\mu_0}:=\left\{\begin{pmatrix}
\phi \\ \psi \end{pmatrix}\,:\, \phi\in L^1_{loc}((-\infty,0]),\,\psi\in C((-\infty,0]),\,\,\phi\ge 0,\,\psi\ge 0,\,\,\left|\left|\begin{pmatrix}
\phi \\ \psi \end{pmatrix} \right|\right|_{\mu_0}<\infty\right\},
\end{equation*}
and it is equipped with the norm $||\cdot||_{\mu_0}$. We refer to $\mathcal{X}_{\mu_0}$ as the space of birth rates and environmental histories. 
\end{definition}
Invoking results from \cite{DG} we have.
\begin{theorem} 
$\left\{\mathcal{T}_{DE}(t)\right\}_{t\ge 0}$ is a strongly continuous semigroup of nonlinear operators on $\mathcal{X}_{\mu_0}$, and the principle of linearised stability holds for this semigroup, 
whenever the nonlinear maps $F_1$ and $F_2$, corresponding to the right-hand side of (\ref{DE-1}) and defined precisely in (\ref{sect6-eq1}) below, are $C^1$.
\end{theorem}
Note that in Appendix \ref{appendix-diffF} we prove differentiability of the maps $F_1,F_2$ under some hypotheses.

\section{The relationship between the two formulations}

In this section we present continuous functions mapping orbits of one formulation to orbits of the other formulation. The relation is not one to one because many initial pairs of histories $(\phi,\psi)\in\mathcal{X}_{\mu_0}$ determine the same orbit in the space of population densities and environmental conditions. 

\begin{lemma}\label{lemmaBoundsXandT}
Suppose that assumptions H2$_g$ and H3$_g$ hold, i.e. 
\begin{equation}\label{growth-assum}
0 < g_{\text{min}}\leq g(x,S)\leq g_{\text{max}},\qquad\forall x\in[x_b,\infty),\quad S\in \mathbb{R}_+.
\end{equation}
Moreover assume that H$_{g_\infty}$ holds. Then for $\theta\leq 0$, there exist constants $c_1$ and $c_2$, such that
\begin{equation}\label{X-upper-bound}
c_1-g_\infty\theta \leq X_\psi(0,\theta,x_b)\leq c_2 - g_\infty \theta,
\end{equation}
and for $x\geq x_b$, we have
\begin{equation}\label{tau-upper-bound}
\frac{c_1-x}{g_\infty}\leq T_\psi(x_b,x,0) \leq \frac{c_2-x}{g_\infty}.
\end{equation}
A possible choice for $c_1$ and $c_2$ is
\[
c_1=\bar{x}-\frac{\bar{x}-x_b}{g_\text{min}}g_\infty\qquad\text{and}\qquad c_2=\bar{x}.
\]  
\end{lemma}
{\bf Proof.}
Assumption \eqref{growth-assum} implies that $X_\psi(0,\theta,x_b)>\bar{x}$ for small enough $\theta$, specifically if $\theta<\tilde\theta=-(\bar{x}-x_b)/g_\text{min}$. This fact together with assumption H$_{g_\infty}$ imply that $X_\psi(0,\theta,x_b)=c-\theta g_\infty$ if $\theta<\tilde{\theta}$, where $c\in [\bar{x}+\tilde{\theta}g_\infty,\bar{x}]$.

Using $X_\psi(0,T_\psi(x_b,x,0),x_b)=x$ in \eqref{X-upper-bound}, the bounds in \eqref{tau-upper-bound} are obtained.

\eofproof

For any weight $\mu_0\in(0,\hat\mu]$ determining the space of histories (associated to the DE formulation) choose $\kappa_0=(\hat\mu-\mu_0)/g_\infty$ as the counterpart weight in the space of densities (associated to the PDE formulation). Notice that for $\mu_0=\hat{\mu}$, we have $\kappa_0=0$. For this special case, the equivalence between the two formulations can be proven without using assumption H$_{g_\infty}$. This assumption is used only to prove the equivalence in the other cases, i.e. when $0<\mu_0<\hat{\mu}$ and $\kappa_0=(\hat{\mu}-\mu_0)/g_\infty>0$. As we already mentioned, in Appendix \ref{appendix-diffF} we find that to prove the differentiability of the delay equation, we need a bound on $\mu_0$, viz.  $\mu_0<\frac{\hat{\mu}}{3}$ (note that this bound is not necessarily sharp).


\begin{theorem}\label{contiuity-L}
Assume H$_s$, H2$_g$, H3$_g$ and H$_{g_\infty}$. Then the map
\begin{equation*}
\begin{array}{rccc}
\mathcal{L}: & \mathcal{X}_{\mu_0} & \longrightarrow & L^1_{\kappa_0}\times\mathbb{R}_+ \\
 & \begin{pmatrix}
\phi \\ \psi \end{pmatrix} & \longmapsto & \begin{pmatrix}
n_0 \\ \psi(0) \end{pmatrix}
\end{array}
\end{equation*}
with $n_0$ given `explicitly' by \eqref{DE-2} is, for $\kappa_0=\frac{\hat{\mu}-\mu_0}{g_\infty}$, well defined (i.e. $\mathcal{L}(\mathcal{X}_{\mu_0})\subset L^1_{\kappa_0}\times\mathbb{R}_+$) and continuous.
\end{theorem}

{\bf Proof.} (of the first half of the statement. The continuity of $\mathcal{L}$ is proven in Appendix \ref{appendix-continuityL}.)
Let us check that $n_0$ defined by \eqref{DE-2} belongs to $L^1_{\kappa_0}$. Indeed
\begin{equation*}
\int_{x_b}^{\infty} n_0 (x) e^{\kappa_0 x}\,\ud x = \int_{x_b}^{\infty} \phi
   (T_\psi(x_b,x,0)) \mathcal{F}_{\psi} (0, T_\psi(x_b,x,0), x_b) (- D_2
   T_\psi(x_b,x,0)) e^{\kappa_0 x}\,\ud x,
\end{equation*}
and making the change of variables $\theta = T_\psi(x_b,x,0)$ (using $X_{\psi} (0, T_\psi(x_b,x,0), x_b) = x$), we obtain
\begin{equation*}
- \int_{T_\psi(x_b,x_b,0)}^{T_\psi(x_b,\infty,0)} \phi (\theta)
   \mathcal{F}_{\psi} (0, \theta, x_b) e^{\kappa_0 X_{\psi} (0, \theta, x_b)}\,\ud \theta =
   \int_{- \infty}^0 \phi (\theta) \mathcal{F}_{\psi} (0, \theta, x_b) e^{\kappa_0
   X_{\psi} (0, \theta, x_b)}\, \ud \theta ,
\end{equation*}
where in the last step we used $T_\psi(x_b,\infty,0)=-\infty$, cf. Lemma \ref{lemmaBoundsXandT}.

Then, applying assumption (\ref{weight-ass1}) and the upper bound of (\ref{X-upper-bound}) one has
\begin{equation*}
\int_{- \infty}^0 \phi (\theta) \mathcal{F}_{\psi} (0, \theta, x_b) e^{\kappa_0
   X_{\psi} (0, \theta, x_b)}\,\ud \theta \leqslant C \int_{- \infty}^0 \phi
   (\theta) e^{\hat{\mu} \theta} e^{\kappa_0 (c_2 - g_{\infty} \theta)}\,\ud \theta
   \leqslant C_1 \int_{- \infty}^0 \phi (\theta) e^{(\hat{\mu} - \kappa_0
   g_{\infty}) \theta}\,\ud \theta
\end{equation*}
so that, using $\kappa_0 g_{\infty}= \hat{\mu} - \mu_0$, we conclude
\begin{equation*}
\int_{x_b}^{\infty} n_0 (x) e^{\kappa_0 x}\,\ud x \leqslant C_1 \int_{-
   \infty}^0 \phi (\theta) e^{\mu_0 \theta}\,\ud \theta \leqslant C_1 \| \phi \|_{1,\mu_0}
   < \infty.
\end{equation*}
\eofproof
   
As an immediate consequence of the constructive definition of both $\mathcal{T}_{PDE}(t)$ and $\mathcal{T}_{DE}(t)$ we have the following result.
\begin{theorem}\label{estimate-th}
\begin{equation}\label{L-th}
\mathcal{T}_{PDE}(t)\,\mathcal{L}=\mathcal{L}\,\mathcal{T}_{DE}(t),\quad t\ge 0.
\end{equation}
\end{theorem}
This result motivates us to introduce an equivalence relation, denoted by 
"$\sim$", on $\mathcal{X}_{\mu_0}$ as follows.

\begin{definition}
We write $\begin{pmatrix}
\phi^1 \\ \psi^1 \end{pmatrix}\sim\begin{pmatrix}
\phi^2 \\ \psi^2 \end{pmatrix}$, and say that these pairs of functions are equivalent, if and only if 
$$\mathcal{L}\begin{pmatrix}
\phi^1 \\ \psi^1 \end{pmatrix}=\mathcal{L}\begin{pmatrix}
\phi^2 \\ \psi^2 \end{pmatrix}.$$
\end{definition}

\begin{theorem}
If $\begin{pmatrix}
\phi^1 \\ \psi^1 \end{pmatrix}\sim\begin{pmatrix}
\phi^2 \\ \psi^2 \end{pmatrix}$ then $\mathcal{T}_{DE}(t)\begin{pmatrix}
\phi^1 \\ \psi^1 \end{pmatrix}\sim\mathcal{T}_{DE}(t)\begin{pmatrix}
\phi^2 \\ \psi^2 \end{pmatrix}$ for $t\ge 0$. Moreover, we have
\begin{equation*}
\left|\left|\mathcal{T}_{DE}(t)\begin{pmatrix}
\phi^1 \\ \psi^1 \end{pmatrix}-\mathcal{T}_{DE}(t)\begin{pmatrix}
\phi^2 \\ \psi^2 \end{pmatrix}\right|\right|\le e^{-\mu_0 t}\left|\left|\begin{pmatrix}
\phi^1 \\ \psi^1 \end{pmatrix}-\begin{pmatrix}
\phi^2 \\ \psi^2 \end{pmatrix}\right|\right|,\quad \forall\,t\ge 0.
\end{equation*}
\end{theorem} 
{\bf Proof.}
If $\begin{pmatrix}
\phi^1 \\ \psi^1 \end{pmatrix}\sim\begin{pmatrix}
\phi^2 \\ \psi^2 \end{pmatrix}$ then $b^1(t)=b^2(t),\,S^1(t)=S^2(t)$ holds for $t>0$ and therefore we have
\begin{align*}
||b^1_t-b^2_t||_{1,\mu_0}=& \int_t^\infty\left|\phi^1(t-a)-\phi^2(t-a)\right|e^{-\mu_0 a}\,\ud a=e^{-\mu_0 t}||\phi^1-\phi^2||_{1,\mu_0}, \\
\left|\left|S^1_t-S^2_t\right|\right|_{\infty,\mu_0}= & \sup_{t\le a<\infty}|\psi^1(t-a)-\psi^2(t-a)|e^{-\mu_0 a}=e^{-\mu_0 t}||\psi^1-\psi^2||_{\infty,\mu_0}.
\end{align*}
\eofproof
\begin{corollary}
If $\begin{pmatrix}
\phi^1 \\ \psi^1 \end{pmatrix}\sim\begin{pmatrix}
\phi^2 \\ \psi^2 \end{pmatrix}$, then $\begin{pmatrix}
\phi^1 \\ \psi^1 \end{pmatrix}$ and $\begin{pmatrix}
\phi^2 \\ \psi^2 \end{pmatrix}$ have the same $\omega$-limit set with respect to $\mathcal{T}_{DE}(t)$.
\end{corollary}
The motivation to introduce the equivalence relation "$\sim$" is that in general $\mathcal{L}$ is many-to-one (we have already observed this phenomenon in case of the distributed states at birth model in \cite{CDF}). But perhaps one can work with an  appropriate representative of each equivalence class; and if, in addition, $\mathcal{L}$ is surjective, then we can define a pseudo-inverse of $\mathcal{L}$, i.e. a map 
\begin{equation*}
\mathcal{L}^{-1}_{ps}\,:\,L^1_{\kappa_0}\times\mathbb{R}_+\to\mathcal{X}_{\mu_0},
\end{equation*}
such that $\mathcal{L}\,\mathcal{L}^{-1}_{ps}$ is the identity on $L^1_{\kappa_0}\times\mathbb{R}_+$.

To show that a function $\mathcal{L}^{-1}_{ps}$ having the properties mentioned above does exist, first notice that most of the many-to-one character of $\mathcal{L}$ seems to be due to the arbitrariness in the history $\psi$ of the resource concentration. In this component, all functions having the same value at $\theta=0$ are equivalent. It somehow seems natural to choose a constant function as a representative, in particular to facilitate the discussion of steady states. So given $S_0\in\mathbb{R}_+$ we choose
\begin{equation}\label{psi-def}
\psi(\theta)= S_0,\quad \theta\le 0,
\end{equation}
and next focus our attention on \eqref{DE-2}, but now we consider $n_0$ as given and $\phi$ as to be determined. According to \eqref{tau-def} and \eqref{tau-eq}, the transformation
\begin{equation}\label{equiv-1}
T_\psi(x_b,x,0)=\theta
\end{equation}
has inverse
\begin{equation}\label{equiv-2}
x=X_\psi(0,\theta,x_b),
\end{equation}
and
\begin{equation}\label{equiv-3}
D_2\,X_\psi(0,T_\psi(x_b,x,0),x_b)D_2\,T_\psi(x_b,x,0)=1.
\end{equation}
These observations allow us to rewrite \eqref{DE-2} as
\begin{equation}\label{equiv-4}
\phi(\theta)=n_0\left(X_\psi(0,\theta,x_b)\right)\frac{1}{\mathcal{F}_\psi(0,\theta,x_b)}\left(-D_2\, X_\psi(0,\theta,x_b)\right),\quad \theta\le 0.
\end{equation}
That the pair $\psi$ and $\phi$ defined by \eqref{psi-def} and \eqref{equiv-4} in terms of $S$ and $n_0$ gives the desired pseudo-inverse is the content of the following result.

\begin{theorem}\label{pseudo-th}Assume H$_s$, H2$_g$, H3$_g$ and H$_{g_\infty}$. Then the map
\begin{equation}\label{equiv-5}
\begin{array}{rccc}
\mathcal{L}^{-1}_{ps}: & L^1_{\kappa_0}\times\mathbb{R}_+ & \longrightarrow & \mathcal{X}_{\mu_0} \\
 & \begin{pmatrix}
n_0 \\ S_0 \end{pmatrix} & \longmapsto & \begin{pmatrix}
\phi \\ \psi \end{pmatrix}
\end{array}
\end{equation}
with first $\psi$ given by \eqref{psi-def} and next $\phi$ given by \eqref{equiv-4}, for $\mu_0=\hat{\mu}-g_\infty\kappa_0$, is well defined (i.e. $\mathcal{L}^{-1}_{ps} (L^1_{\kappa_0}\times\mathbb{R}_+) \subset \mathcal{X}_{\mu_0}$), continuous and satisfies
\begin{equation}\label{L-id}
\mathcal{L}\,\mathcal{L}^{-1}_{ps}=\mathcal{I}.
\end{equation}
\end{theorem}

{\bf Proof.} (The proof of the continuity of $\mathcal{L}_\text{ps}^{-1}$ is given in Appendix \ref{appendix-continuityL}.)
If $\mathcal{L}^{-1}_{ps}$ is well defined then it satisfies (\ref{L-id}) 
by construction. Therefore, it is enough to show that $\phi$ defined by \eqref{equiv-4} (with $\psi(\cdot)\equiv S_0$) satisfies
\begin{equation*} 
\int_{- \infty}^0 \phi (\theta) e^{\mu_0 \theta}\,\ud \theta < \infty,
\end{equation*}
so that the pair $(\phi,\psi)\in\mathcal{X}_{\mu_0}$.
Indeed
\begin{equation*}
\int_{- \infty}^0 \phi (\theta) e^{\mu_0 \theta} \,\ud \theta = \int_{-
   \infty}^0 n_0 (X_{\psi} (0, \theta, x_b)) \frac{1}{F_{\psi} (0, \theta,
   x_b)} (- D_2 X_{\psi} (0, \theta, x_b)) e^{\mu_0 \theta} \,\ud \theta
\end{equation*}
and making the change of variables $X_{\psi} (0, \theta, x_b) = x$ (using $T_\psi(x_b,X_\psi(0,\theta,x_b),0) = \theta$), it follows that
\begin{equation*}
- \int_{X_{\psi} (0, - \infty, x_b)}^{X_{\psi} (0, 0, x_b)} n_0 (x)
   \frac{ e^{\mu_0 T_\psi(x_b,x,0)}}{F_{\psi} (0, T_\psi(x_b,x,0), x_b)} \,\ud x = \int_{x_b}^{\infty} n_0 (x) \frac{e^{\mu_0 T_\psi(x_b,x,0)}}{F_{\psi} (0,
   T_\psi(x_b,x,0), x_b)}  \,\ud x .
\end{equation*}
Then, applying the lower bound in (\ref{weight-ass1}) and the lower bound 
of (\ref{tau-upper-bound}) one has
\begin{equation*}
\int_{x_b}^{\infty} n_0 (x) \frac{e^{\mu_0 T_\psi(x_b,x,0)}}{F_{\psi} (0, T_\psi(x_b,x,0),x_b)}  \,\ud x \leqslant \frac{1}{c}
     \int_{x_b}^{\infty} n_0 (x) e^{-(\hat{\mu} - \mu_0)T_\psi(x_b,x,0)}  \,\ud x \leqslant C_2 \int_{x_b}^{\infty} n_0 (x) e^{(\hat{\mu} - \mu_0)\frac{x}{g_{\infty}}} \,\ud x,
\end{equation*}
and using $\kappa_0 = \frac{\hat{\mu} - \mu_0}{g_{\infty}}$ we conclude that
\begin{equation*} 
\int_{- \infty}^0 \phi (\theta) e^{\mu_0 \theta} \,\ud \theta \leqslant C_2\int_{x_b}^{\infty} n_0 (x) e^{\kappa_0 x} \,\ud x \leqslant C_2 \|n_0 \|_{L^1_{\kappa_0}} < \infty.
\end{equation*}
\eofproof

Note that the continuity of the `reformulation' maps $\mathcal{L}$ and $\mathcal{L}^{-1}_{ps}$ is needed when we want to transfer stability assertions from one formulation to the other .

\begin{remark} As we already said, assumption H$_{g_\infty}$ is not needed to prove that $\mathcal{L}$ and $\mathcal{L}^{-1}_{ps}$ are well defined and continuous in the case $\mu_0=\hat{\mu}$ and $\kappa_0=0$. Whether this assumption can also be discarded when $\mu_0\in(0,\hat{\mu})$ and $\kappa_0 = (\hat{\mu}-\mu_0)\tilde{g}^{-1}$ for a suitable $\tilde{g}\in[g_{\min} ,g_{\max}]$,  is an open question.\\
\end{remark}

To apply the pseudo-inverse makes sense if we are dealing with an arbitrary element of $L^1_{\kappa_0}([x_b,\infty);\mathbb{R}_+)\times\mathbb{R}_+$. But for 
points on an orbit of $\mathcal{T}_{PDE}(\cdot)$, there exist function $\tau\mapsto (b(\tau),S(\tau))$ for a certain interval of values of $\tau$, 
and we should use this information. We now show that full orbits, i.e. orbits that go back in time to $-\infty$, of $\mathcal{T}_{PDE}(\cdot)$ and 
$\mathcal{T}_{DE}(\cdot)$ are in one-to-one correspondence.

\begin{theorem}\label{equiv-th}
(i) Let $t\mapsto \begin{pmatrix}
b(t) \\ S(t)
\end{pmatrix}$, from $ \mathbb{R}$ to $ \mathbb{R}_+ \times \mathbb{R}_+ $, 
be such that 
\begin{equation}\label{equiv-6}
\begin{pmatrix}
b_t \\ S_t
\end{pmatrix}=\mathcal{T}_{DE}(t-s)\begin{pmatrix}
b_s \\ S_s
\end{pmatrix},\quad\forall\, t,s \quad\ \text{with} \quad t>s.
\end{equation}
Define
\begin{equation}\label{equiv-7}
\begin{pmatrix}
n(t,\cdot) \\ S(t)
\end{pmatrix}=\mathcal{L}\begin{pmatrix}
b_t \\ S_t
\end{pmatrix},
\end{equation}
then 
\begin{equation}\label{equiv-8}
\begin{pmatrix}
n(t,\cdot) \\ S(t)
\end{pmatrix}=\mathcal{T}_{PDE}(t-s)\begin{pmatrix}
n(s,\cdot) \\ S(s)
\end{pmatrix},\quad \forall\,t,s\quad \text{with}\quad t>s.
\end{equation}

(ii) Let $t\mapsto\begin{pmatrix}
n(t,\cdot) \\ S(t)
\end{pmatrix}$, from $\mathbb{R}$ to $
L^1_{\kappa_0}([x_b,\infty);\mathbb{R}_+)\times \mathbb{R}_+$, be such that \eqref{equiv-8} holds. Define
\begin{equation}\label{equiv-9}
b(t)=\int_{x_b}^\infty\beta(\xi,S(t))n(t,\xi)\,\ud\xi.
\end{equation}
Then \eqref{equiv-6} holds.
\end{theorem}
{\bf Proof.} (i) Apply $\mathcal{L}$ to the identity \eqref{equiv-6} and use \eqref{L-th} and \eqref{equiv-7}, then \eqref{equiv-8} follows.

(ii) The key point here is that $S(t)$ is known for $-\infty<t<\infty$, so at any time we can, for any size, determine the time of birth, and in \eqref{char-sol} we can restrict our attention to the formula that expresses $n(t,x)$ in terms of $b$. Inserting this expression into the right hand side of \eqref{equiv-9} we obtain
\begin{equation*}
b(t)=\int_{x_b}^\infty\beta(\xi,S(t))b(T_S(x_b,\xi,t))\mathcal{F}_S(t,T_S(x_b,\xi,t),x_b)(-D_2 T_S(x_b,\xi,t))\,\ud\xi.
\end{equation*}
We then replace the integration variable $\xi$ by the integration variable $a$ defined by
\begin{equation*}
T_S(x_b,\xi,t)=t-a,
\end{equation*}
while noting that this definition entails the identity
\begin{equation*}
X_S(t,t-a,x_b)=\xi.
\end{equation*}
We then have
\begin{equation*}
b(t)=\int_0^\infty\beta(X_S(t,t-a,x_b)S(t))b(t-a)\mathcal{F}_S(t,t-a,x_b)\,\ud a,
\end{equation*}
where we have used that
\begin{equation*}
-D_2\,T_S(x_b,X_S(t,t-a,x_b),t)(-D_2\,X_S(t,t-a,x_b))=1,
\end{equation*}
since
\begin{equation*}
T_S(x_b,X_S(t,t-a,x_b),t)=t-a.
\end{equation*}
Thus we obtained the first equation of \eqref{DE-1}. The second equation of \eqref{DE-1} is derived in exactly the same manner, one just has to replace $\beta$ by $\gamma$.

\eofproof

Since the $\omega$-limit set consists of full orbits, we can use Theorem \ref{equiv-th} to switch back and forth between the PDE and DE formulation, when dealing with elements of an $\omega$-limit set.

To transfer information about stability obtained in the delay formulation 
to the PDE formulation, we can use the identity
\begin{equation}\label{equiv-10}
\mathcal{T}_{PDE}(t)=\mathcal{L}\,\mathcal{T}_{DE}(t)\,\mathcal{L}^{-1}_{ps},
\end{equation}
that follows directly from \eqref{L-th} and \eqref{L-id}. But how about the other direction, i.e. to transfer information concerning stability from the PDE formulation to the DE formulation. If we apply $\mathcal{T}_{PDE}(t)$ to $\begin{pmatrix}
n_0 \\ S_0
\end{pmatrix}$, we construct a solution as described in Section 2, so there exist functions $\tau\mapsto b(\tau)$ and $\tau\mapsto S(\tau)$ defined on $[0,t]$ that we can use. Therefore we define
\begin{equation}\label{equiv-11}
\left(\mathcal{L}^{-1}_{t,ps}\,\mathcal{T}_{PDE}(t)\begin{pmatrix}
n_0 \\ S_0
\end{pmatrix}\right)(\theta)=\left\{\begin{aligned} &\begin{pmatrix}
b(t+\theta) \\ S(t+\theta)
\end{pmatrix},& \quad \text{if}\quad t+\theta\ge 0, \\ 
&\mathcal{L}^{-1}_{0,ps}\begin{pmatrix}
n_0 \\ S_0
\end{pmatrix}(t+\theta),& \quad\text{if}\quad t+\theta\le 0, \end{aligned} 
\right.
\end{equation}
where we have written $\mathcal{L}^{-1}_{0,ps}$ to denote $\mathcal{L}^{-1}_{ps}$ defined by \eqref{equiv-5}.

Exactly as in the proof of Theorem \ref{estimate-th} this yields the estimate

\begin{equation*}
\left|\left| \mathcal{L}^{-1}_{t,ps}\,\mathcal{T}_{PDE}(t) \mathcal{L}\begin{pmatrix}
\phi \\ \psi
\end{pmatrix}-\mathcal{T}_{DE}(t)\begin{pmatrix}
\phi \\ \psi
\end{pmatrix}\right|\right|\le e^{-\mu_0 t}\left|\left|\begin{pmatrix}
\phi \\ \psi
\end{pmatrix}-\mathcal{L}^{-1}_{0,ps}\mathcal{L}\begin{pmatrix}
\phi \\ \psi
\end{pmatrix}\right|\right| .
\end{equation*}

In the next section we show that $\mathcal{L}$ sends steady states $(\phi_\ast,\psi_\ast)$ in $\mathcal{X}_{\mu_0}$ to steady states $(n_\ast,S_\ast)$ in $L^1_{\kappa_0}\times\mathbb{R}_+$ and that, reciprocally, $\mathcal{L}_{ps}^{-1}$ sends steady states $(n_\ast,S_\ast)$ in $L^1_{\kappa_0}\times\mathbb{R}_+$ to the steady states $(\phi_\ast,\psi_\ast)$ in $\mathcal{X}_{\mu_0}$, whose image by $\mathcal{L}$ is $(n_\ast,S_\ast)$. Moreover, we also show that the intertwined pairs of steady states $(\phi_\ast,\psi_\ast)$ and $(n_\ast,S_\ast)$ are either both stable or both unstable.

\section{Steady states and stability}

If the resource concentration has a constant value, say $S_*$, then 

\begin{equation*}
X_{S_*}(t,s,x_b)=\overline{\xi}(t-s),
\end{equation*}
where $\overline{\xi}$ is defined as the solution of the ODE
\begin{equation} \label{sect5-eq1}
\dot{\overline{\xi}}=g\left(\overline{\xi},S_*\right),\quad \overline{\xi}(0)=x_b, 
\end{equation}
(note that $\overline{\xi}$ depends on $S_*$ even though this is not expressed in the notation; also note that $S_*$ is sometimes used to denote the constant function taking the value $S_*$). Accordingly, the renewal equation for $b$ as given in \eqref{DE-1} at a steady state reduces to 
\begin{equation} \label{sect5-eq2}
b(t)=\int_0^\infty\beta\left(\overline{\xi}(a),S_*\right)\exp\left\{-\int_0^a\mu\left(\overline{\xi}(\tau),S_*\right)\,\ud \tau\right\} b(t-a)\,\ud a.
\end{equation}
This linear equation has constant solutions if and only if 
\begin{equation}\label{sect5-eq3}
\mathcal{R}\left(S_*\right)=1,
\end{equation}
where (using the transformations $\bar{\xi}(a)=x$ and $\bar{\xi}(\tau)=y$ to obtain the second expression)
\begin{equation}\label{sect5-eq4}
\begin{aligned}
\mathcal{R}\left(S_*\right):= & \int_0^\infty\beta\left(\overline{\xi}(a),S_*\right)\exp\left\{-\int_0^a\mu\left(\overline{\xi}(\tau),S_*\right)\,\ud \tau\right\}\ud a \\
= & \int_{x_b}^\infty \frac{\beta(x,S_*)}{g(x,S_*)}\exp\left\{-\int_{x_b}^x\frac{\mu\left(y,S_*\right)}{g(y,S_*)}\ud y\right\}\ud x;
\end{aligned}
\end{equation}
is the expected number of offspring produced by a newborn individual in the environment characterised by constant resource availability $S_*$. 

Note that \eqref{sect5-eq3} is one equation in one unknown. For certain classes of reasonable functions describing growth, survival and reproduction, $\mathcal{R}$ is monotone increasing and $\mathcal{R}(0)=0$. If this is the case then obviously  \eqref{sect5-eq3} has a unique solution whenever we can find a feasible value $S_0$ such that 
\begin{equation} \label{sect5-eq5}
\mathcal{R}(S_0)>1.
\end{equation} 
A natural candidate for such an $S_0$ is a stable steady state of the ODE
\begin{equation} \label{sect5-eq6}
\frac{\ud S}{\ud t}(t)=f(S(t)).
\end{equation}
Inequality \eqref{sect5-eq5} then simply means that the consumer population starts to grow when introduced in a `virgin' environment (or, in other 
words, that the trivial steady state $(b,S)=(0,S_0)$ is unstable).

If \eqref{sect5-eq3} holds, then every constant function $b(t)\equiv b_*$ 
satisfies \eqref{sect5-eq2}. The `right' $b_*$ is determined from the requirement that population level consumption matches resource production, that is
\begin{equation} \label{sect5-eq7}
f(S_*)=b_*\int_0^\infty\gamma\left(\overline{\xi}(a),S_*\right) \exp\left\{-\int_0^a\mu\left(\overline{\xi}(\tau),S_*\right)\,\ud \tau\right\}\ud a.
\end{equation}  
Indeed, \eqref{sect5-eq7} guarantees that $\frac{\ud S}{\ud t}$ in \eqref{DE-1} equals zero when $b(t)\equiv b_*,\, S(t)\equiv S_*$ for all $t$. Also note that 
\begin{equation*}
\int_0^\infty\gamma\left(\overline{\xi}(a),S_*\right) \exp\left\{-\int_0^a\mu\left(\overline{\xi}(\tau),S_*\right)\,\ud \tau\right\}\ud a=
\int_{x_b}^\infty \frac{\gamma(x,S_*)}{g(x,S_*)}\exp\left\{-\int_{x_b}^x\frac{\mu\left(y,S_*\right)}{g(y,S_*)}\ud y\right\}\,\ud x,
\end{equation*}
is the expected lifetime consumption of resources of a newborn individual, given constant resource concentration $S_*$.

For constant resource concentration $S_*$, the time $T_{S_*}(x_b,x,t)$ of birth of an individual having size $x$ at time $t$, as introduced in 
\eqref{tau-def}-\eqref{tau-eq}, is given by
\begin{equation} \label{sect5-eq8}
T_{S_*}(x_b,x,t)=t-\int_{x_b}^x\frac{\ud y}{g(y,S_*)}.
\end{equation}
At a steady state the second expression in \eqref{char-sol} yields 
\begin{equation} \label{sect5-eq9}
n_*(x)=\frac{b_*}{g(x,S_*)}\exp\left\{-\int_{x_b}^x\frac{\mu\left(y,S_*\right)}{g(y,S_*)}\,\ud y\right\},
\end{equation}
which, together with 
\begin{equation} \label{sect5-eq10}
b_*=\frac{f(S_*)}{\displaystyle\int_{x_b}^\infty \frac{\gamma(x,S_*)}{g(x,S_*)}\exp\left\{-\int_{x_b}^x\frac{\mu\left(y,S_*\right)}{g(y,S_*)}\,\ud y\right\}\,\ud x}
\end{equation}
(see \eqref{sect5-eq7}), yields an explicit formula for the steady size distribution $n_*$, once $S_*$ is determined from \eqref{sect5-eq3}.

Next we show that the steady states of the two formulations are in a one to one correspondence given by the operators $\mathcal{L}$ and $\mathcal{L}_{ps}^{-1}$, and that corresponding steady states share the stability character (in the sense that the steady states are either both stable or both unstable).

Let $S_\ast$ satisfy \eqref{sect5-eq3}. Define $b_\ast$ by \eqref{sect5-eq10} and $n_\ast$ by \eqref{sect5-eq9}.

Then
\begin{enumerate}
\item[i.] $(n_\ast,S_\ast)$ is a steady state of $T_{PDE}(t)$,
\item[ii.] $(b_\ast,S_\ast)$ is a steady state of $T_{DE}(t)$ \\
(note that in i. $S_\ast$ is a scalar, while in ii. it is a constant function with value $S_\ast$),
\item[iii.] $\mathcal{L}(b_\ast,S_\ast)=(n_\ast,S_\ast)$,
\item[iv.] $\mathcal{L}_{ps}^{-1}(n_\ast,S_\ast)=(b_\ast,S_\ast)$.
\end{enumerate}

\begin{theorem}\label{stab-eq-th} Assume H1$_\beta$ and H2$_\beta$. The equilibrium $(n_\ast,S_\ast)$ is stable with respect to $T_{PDE}$ if and only if the equilibrium $(b_\ast,S_\ast)$ is stable with respect to $T_{DE}$.
\end{theorem}

{\bf Proof.} We first show the relatively easy ``if'' part. So assume that $(b_\ast,S_\ast)$ is stable, i.e., assume that $\forall \varepsilon>0$, $\exists \,\delta_2=\delta_2(\varepsilon)$, such that
\begin{equation*}
\|(\phi,\psi)-(b_\ast,S_\ast)\|<\delta_2\quad\Rightarrow \quad\forall t\geq 0,\;\; \|T_{DE}(t)(\phi,\psi)-(b_\ast,S_\ast)\|<\varepsilon.
\end{equation*}
The continuity of $\mathcal{L}_{ps}^{-1}$ guarantees that $\forall \varepsilon>0$, $\exists \,\delta_1=\delta_1(\varepsilon)$, such that
\begin{equation*}
\|(n_0,\psi_0)-(n_\ast,S_\ast)\|<\delta_1\quad\Rightarrow \quad \|\mathcal{L}_{ps}^{-1}(n_0,\psi_0)-(b_\ast,S_\ast)\|<\varepsilon,
\end{equation*}
and the continuity of $\mathcal{L}$ guarantees that $\forall \varepsilon>0$, $\exists \,\delta_3=\delta_3(\varepsilon)$ such that
\begin{equation*}
\|(\phi,\psi)-(b_\ast,S_\ast)\|<\delta_3\quad\Rightarrow \quad \|\mathcal{L}(\phi,\psi)-(n_\ast,S_\ast)\|<\varepsilon.
\end{equation*}

Now choose $\delta_4(\varepsilon)=\delta_1(\delta_2(\delta_3(\varepsilon)))$, then
\begin{equation*}
\begin{aligned}
&\quad \|(n_0,\psi_0)-(n_\ast,S_\ast)\|<\delta_4(\varepsilon) \\
\Rightarrow&\quad \|\mathcal{L}_{ps}^{-1}(n_0,\psi_0)-(b_\ast,S_\ast)\|<\delta_2(\delta_3(\varepsilon)) \\
\Rightarrow & \quad t\geq 0,\;\; \|T_{DE}(t)\mathcal{L}_{ps}^{-1}(n_0,\psi_0)-(b_\ast,S_\ast)\|<\delta_3(\varepsilon)\\
\Rightarrow & \quad t\geq 0,\;\; \|\mathcal{L}T_{DE}(t)\mathcal{L}_{ps}^{-1}(n_0,\psi_0)-(n_\ast,S_\ast)\|< \varepsilon,
\end{aligned}
\end{equation*}
and the use of \eqref{equiv-10} completes the proof.

Next we concentrate on the ``only if'' part. Assume that $(n_\ast,S_\ast)$ is stable: $\forall \varepsilon>0$, $\exists \,\delta_5=\delta_5(\varepsilon)$, such that
\begin{equation*}
\|(n_0,\psi_0)-(n_\ast,S_\ast)\|<\delta_5\quad\Longrightarrow \quad \forall t\geq 0,\;\; \|T_{PDE}(t)(n_0,\psi_0)-(n_\ast,S_\ast)\|<\varepsilon.
\end{equation*}
Let $\delta_3=\delta_3(\varepsilon)$ be as above, i.e. a characterisation of the continuity of $\mathcal{L}$.

For given $(\phi,\psi)\in \mathcal{X}_{\mu_0}$ let $b=b(\phi,\psi)$ be the birth rate. Let $S=S(\phi,\psi)$ be the resource concentration.


We shall use
\begin{equation*}
\eqref{birth-rate} \hspace{26mm} b(t)=\int_{x_b}^\infty \beta(\xi,S(t))n(t,\xi)\,\ud\xi,
\end{equation*}
and its steady state version 
\begin{equation*}
\eqref{sect5-eq9},\eqref{sect5-eq4} \hspace{20mm} b_*=\int_{x_b}^\infty \beta(\xi,S_*)n_*(\xi)\,\ud\xi, 
\end{equation*}
and
\begin{itemize}
\item[H2$_\beta$] $\beta(x,S) \leq B_0$,
\item[H1$_\beta$] $ |\beta(x_1,S_1)-\beta(x_2,S_2)|\leq B_1 |x_1-x_2|+B_2|S_1-S_2|$.
\end{itemize}

\begin{lemma}\label{onlyif-lemma} Assume H1$_\beta$ and H2$_\beta$. Let $(n_0,S_0)=\mathcal{L}(\phi,\psi)$. If for all $t\geq 0$
\[
\|T_{PDE}(t)(n_0,S_0)-(n_\ast,S_\ast)\|<\varepsilon,
\]
then
\[
|b(t)-b_\ast|<\varepsilon\left(B_0+B_2\int_{x_b}^\infty n_{\ast}(\xi)\,\ud\xi\right).
\]
\end{lemma}
\textbf{Proof.}
Notice that
\[
\|T_{PDE}(t)(n_0,S_0)-(n_\ast,S_\ast)\|=\int_{x_b}^\infty |n(t,\xi)-n_\ast(\xi)|e^{\kappa_0 \xi}\,\ud\xi+|S(t)-S_\ast|.
\]
Therefore the assumption guarantees that both terms at the right hand side are bounded by $\varepsilon$. Next note that (using $\kappa_0\geq 0)$
\[
\begin{aligned}
|b(t)-b_\ast|&=\left|\int_{x_b}^\infty \beta(\xi,S(t))n(t,\xi)-\beta(\xi,S_\ast)n_\ast(\xi)\,\ud\xi\right| \\
& \leq \int_{x_b}^\infty \beta(\xi,S(t))|n(t,\xi)-n_\ast(\xi)|\,\ud\xi+\int_{x_b}^\infty |\beta(\xi,S(t))-\beta(\xi,S_\ast)|n_\ast(\xi)\,\ud\xi \\
& \leq B_0\int_{x_b}^\infty |n(t,\xi)-n_\ast(\xi)|e^{\kappa_0\xi}\,\ud\xi+B_2|S(t)-S_\ast|\int_{x_b}^\infty n_\ast(\xi)\,\ud\xi\\
& < \varepsilon\left(B_0+B_2\int_{x_b}^\infty n_{\ast}(\xi)\,\ud\xi\right). 
\end{aligned}
\]
\eofproof

Next note that
\begin{equation}\label{Tde-bound}
\|T_{DE}(t)(\phi,\psi)-(b_\ast,S_\ast)\|=(1)+(2)+\sup\{(3),(4)\}
\end{equation}
with
\begin{equation}\label{1234terms}
\begin{aligned}
&(1)=\int_0^t|b(t-a)-b_\ast|e^{-\mu_0 a}\,\ud a, \\
&(2)= \int_t^\infty|\phi(t-a)-b_\ast|e^{-\mu_0 a}\,\ud a, \\
&(3)=\sup_{0\leq a\leq t} |S(t-a)-S_\ast|e^{-\mu_0 a},\\
&(4)=\sup_{ a> t} |\psi(t-a)-S_\ast|e^{-\mu_0 a},
\end{aligned}
\end{equation}
and that
\begin{equation*}
(2)+(4)=e^{-\mu_0 t}\|(\phi,\psi)-(b_\ast,S_\ast)\|\leq \|(\phi,\psi)-(b_\ast,S_\ast)\|,\quad\forall t\geq 0.
\end{equation*}

%

Choose $\delta_6(\varepsilon)=\min\left\{\frac{\varepsilon}{2},\delta_3\left(\delta_5\left(\frac{\varepsilon}{2}Q\right)\right)\right\}$ with
\begin{equation*}
Q:=\left(1+\frac{B_0+B_2\int_{x_b}^\infty n_\ast(\xi)d\xi}{\mu_0}\right)^{-1}.
\end{equation*}
Then, setting $(n_0,S_0)=\mathcal{L}(\phi,\psi)$, it follows that
\begin{equation*}
\begin{aligned}
&\quad \|(\phi,\psi)-(b_\ast,S_\ast)\|<\delta_3\left(\delta_5\left(\frac{\varepsilon}{2}Q\right)\right) \\
\Rightarrow&\quad \|(n_0,S_0)-(n_\ast,S_\ast)\|<\delta_5\left(\frac{\varepsilon}{2}Q\right) \\
\Rightarrow & \quad \forall\,t\geq 0,\;\; \|T_{PDE}(t)(n_0,S_0)-(n_\ast,S_\ast)\|<\frac{\varepsilon}{2}\,Q\\
\overset{\ast}{\Rightarrow} & \quad \forall\,t\geq 0,\;\; |b(t)-b_\ast|<\frac{\varepsilon}{2}\,Q\left(B_0+B_2\int_{x_b}^\infty n_{\ast}(\xi)\,\ud\xi\right),
\end{aligned}
\end{equation*}
where in $\overset{\ast}{\Rightarrow}$ we use Lemma \ref{onlyif-lemma} (with $\frac{\varepsilon}{2}\,Q$ instead of $\varepsilon$). As observed in the Proof of Lemma \ref{onlyif-lemma}, we also have
\begin{equation*}
|S(t)-S_\ast|\leq \frac{\varepsilon}{2}Q.
\end{equation*}
Hence
\begin{equation*}
(1)+(3)<\frac{\varepsilon}{2}\,Q\left(B_0+B_2\int_{x_b}^\infty n_{\ast}(\xi)\,\ud\xi\right)\frac{1}{\mu_0}+\frac{\varepsilon}{2}\,Q=\frac{\varepsilon}{2}.
\end{equation*}
On the other hand
\begin{equation*}
(2)+(4)\leq \|(\phi,\psi)-(b_\ast,S_\ast)\|<\frac{\varepsilon}{2}, \quad\forall t\geq 0,
\end{equation*}
so that, using \eqref{Tde-bound}, if $\|(\phi,\psi)-(b_\ast,S_\ast)\|<\delta_6(\varepsilon)$, then
\begin{equation*}
\forall t\geq 0,\;\;\|T_{DE}(t)(\phi,\psi)-(b_\ast,S_\ast)\|=(1)+(2)+\sup\{(3),(4)\}\leq (1)+(2)+(3)+(4)<\varepsilon.
\end{equation*}
\eofproof

\begin{theorem}
Assume H1$_\beta$ and H2$_\beta$. The equilibrium $(n_\ast,S_\ast)$ is asymptotically stable with respect to $T_{PDE}$ if and only if the equilibrium $(b_\ast,S_\ast)$ is asymptotically stable with respect to $T_{DE}$.
\end{theorem}

\textbf{Proof.} To prove the if part, first use that the asymptotic stability of $(b_\ast,S_\ast)$ by $T_{DE}$ implies that
there exists a ball $B1$ (in $\mathcal{X}_{\mu_0}$) centered at $(b_\ast,S_\ast)$ for which $T_{DE}(t)(b,S)\rightarrow (b_\ast,S_\ast)$ as $t\rightarrow \infty$ for all $(b,S)$ in $B1$. Then take a second ball $B2$ (now in $L_{\kappa_0}^1\times \mathbb{R}$)
centered at $(n_\ast,S_\ast)$ small enough such that $\mathcal{L}_{ps}^{-1}(B2) \subset B1$, which
is possible because of the continuity of $\mathcal{L}_{ps}^{-1}$ and because $\mathcal{L}_{ps}^{-1}(n_\ast,S_\ast)=(b_\ast,S_\ast)$. Then it follows that for all $(n,S)\in B2$, $T_{PDE}(t)(n,S)\rightarrow (n_\ast,S_\ast)$ as $t\rightarrow \infty$. Indeed, since by construction $\mathcal{L}_{ps}^{-1}(n,S) \in B1$, then
\[
\lim_{t\rightarrow \infty} T_{PDE}(t)(n,S)  = \lim_{t\rightarrow \infty} \mathcal{L}( T_{DE}(t)\mathcal{L}_{ps}^{-1}(n,S) ) ) = \mathcal{L} \left(  \lim_{t\rightarrow \infty}
T_{DE}(t)L_{ps}^{-1}(n,S) \right) =  \mathcal{L}(b_\ast,S_\ast) = (n_\ast,S_\ast).
\]

To prove the only if part take a small enough ball $B$ centred at $(b_\ast,S_\ast)$ so that the image of this ball by $\mathcal{L}$ is contained in the basin of attraction of $(n_\ast,S_\ast)$. Then, for each pair of histories
$(\phi,\psi)\in B$, define $(n_0,S_0)=\mathcal{L}(\phi,\psi)$. Since the orbit of $(n_0,S_0)$ tends to $(n_\ast,S_\ast)$ by construction, the function
\[ 
\varepsilon(t):= \|T_{PDE}(t)(n_0,S_0)-(n^\ast,S^\ast)\|
\]
tends to zero as $t\rightarrow 0$. Now, using the same reasoning as in Lemma \ref{onlyif-lemma}, it follows that
\[
|S(t)-S^\ast|\leq \varepsilon(t)\quad\text{\and}\quad |b(t)-b^\ast|\leq \varepsilon(t)\left(B_0+B_2\int_{x_b}^\infty n_{\ast}(\xi)\,\ud\xi\right),
\]
which implies the only if part since 
\[
\|T_{DE}(t)(\phi,\psi)-(b*,S*)\| = (1)+(2)+\sup\{(3)+(4)\} \leq (1)+(3)+(2)+(4),
\]
with (1),(2),(3) and (4) defined in (\ref{1234terms}), and (1)+(3) tends to 0 because $\varepsilon(t)$ tends to 0 and (2)+(4) decays to
0 exponentially at a rate $\mu_0$.

\eofproof

The standard procedure to determine the stability of steady states involves linearisation of the solution operators. When the equations can be linearised, the linearised solution operators are obtained as the solution operators of the linearised equations. In the PDE formulation, formal linearisation of the equations is no problem at all and the corresponding eigenvalue problem does indeed lead to the correct characteristic equation \eqref{sect6-eq9} below. But note that when linearising the growth term, we formally differentiate an unbounded operator. To prove the Principle of Linearised Stability along this route is impossible, for the simple reason that the solution operators are, in fact, not differentiable.
 Indeed, as already noted at the end of Section \ref{sectPDEformulation}, the first part of \eqref{char-sol} shows that solving the problem involves shifting the initial function $n_0$ over an $S$-dependent distance, and when $n_0$ is not absolutely continuous this operation does not depend differentiably on $S$ (this is, we recall, a manifestation of the smoothness problem created by state-dependent delay).
   The non-differentiability is transient : the $x$-domain, to which the first formula of \eqref{char-sol} applies, shifts towards infinity when time proceeds and the norm of the corresponding part of the solution decays exponentially under natural conditions on the per capita death rate. The persistent behaviour is, hopefully, described by differentiable operators. 
 It is exactly when dealing with linearised stability that the delay formulation offers an advantage: certain assumptions on the model ingredients guarantee that the solution operators are differentiable in the delay setting for suitably chosen $\mu_0$. We show this in Appendix C.

Here we formally derive the linearisation of the DE formulation \eqref{DE-1}. First, in order to be able to apply the results in Section 5 of \cite{DG}, we need to write \eqref{DE-1} in the form
\begin{equation} \label{sect6-eq1}
\begin{aligned}
b(t)= & F_1(b_t,S_t), \\
\frac{\ud S}{\ud t}(t)= & F_2(b_t,S_t),
\end{aligned}
\end{equation}
and check that $F_1,\, F_2\,:\,\mathcal{X}\to\mathbb{R}$ are $C^1$ under appropriate conditions on $\beta, \mu, g, \gamma$ and $f$. To this end, we first observe that from \eqref{X-s-def-2} and \eqref{F-s-def-2} it follows that (recall the notation \eqref{DE-5}) 
\begin{equation}\label{sect6-eq2}
\begin{aligned}
X_S(t,t-a,x_b)= & X_{S_t}(0,-a,x_b),\\
\mathcal{F}_S(t,t-a,x_b)= & \mathcal{F}_{S_t}(0,-a,x_b).
\end{aligned}
\end{equation}
So if we define $F_1$ and $F_2$ as
\begin{equation}\label{sect6-eq3}
\begin{aligned}
F_1(\phi,\psi)= & \int_0^\infty\beta\left(X_\psi(0,-a,x_b),\psi(0)\right)\mathcal{F}_\psi(0,-a,x_b)\phi(-a)\,\ud a,\\
F_2(\phi,\psi)= & f(\psi(0))-\int_0^\infty\gamma\left(X_\psi(0,-a,x_b),\psi(0)\right)\mathcal{F}_\psi(0,-a,x_b)\phi(-a)\,\ud a,
\end{aligned}
\end{equation}
then \eqref{DE-1} does indeed correspond to \eqref{sect6-eq1}. Note that $F_1$ and $F_2$ are well-defined if 
\begin{equation}\label{F-cond-bound}
\sup_{a\ge 0}\left\{ e^{\mu_0a}\beta\left(X_\psi(0,-a,x_b),\psi(0)\right)\mathcal{F}_\psi(0,-a,x_b)\right\}<\infty,\quad \sup_{a\ge 0}\left\{ e^{\mu_0a}\gamma\left(X_\psi(0,-a,x_b),\psi(0)\right)\mathcal{F}_\psi(0,-a,x_b)\right\}<\infty.
\end{equation}

In Appendix \ref{appendix-diffF} Theorem \ref{theoremC12} we prove that $F_1,\, F_2$ map $\mathcal{X}_{\mu_0}$ into $\mathbb{R}$ and are indeed of class $C^1$. It follows that under these assumptions the Principle of Linearised Stability holds.

Let $(b_*,S_*)$ be a non-trivial steady state, i.e. assume that 
\begin{equation} \label{sect6-eq4}
\begin{aligned}
b_*= & F_1(b_*,S_*), \\
0= & F_2(b_*,S_*).
\end{aligned}
\end{equation}
If we insert
\begin{equation} \label{sect6-eq5}
\begin{aligned}
b(t)= & b_*+\varepsilon y(t), \\
S(t)= & S_*+\varepsilon z(t)
\end{aligned}
\end{equation}
into \eqref{sect6-eq1}, divide by $\varepsilon$ and let $\varepsilon\to 0$, we obtain the linearised system
\begin{equation} \label{sect6-eq6}
\begin{aligned}
y(t)= & D_1\,F_1(b_*,S_*)\, y_t+D_2\, F_1(b_*,S_*)\, z_t, \\
\frac{\ud z}{\ud t}(t)= & D_1\, F_2(b_*,S_*)\, y_t+D_2\, F_2(b_*,S_*)\, 
z_t.
\end{aligned}
\end{equation}
Note that 
\begin{equation}\label{sect6-eq7}
D_1\, F_i(b_*,S_*)\, y_t=F_i(y_t,S_*),\quad i=1,2,
\end{equation}
since $F_1$ and $F_2$ are linear in the first component. We will derive representations for $D_2\, F_1(b_*,S_*)$ and $D_2\, F_2(b_*,S_*)$ later on.  The linear system \eqref{sect6-eq6} admits a solution of the form
\begin{equation}\label{sect6-eq8}
\begin{pmatrix}
y(t) \\ 
z(t)
\end{pmatrix}= e^{\lambda t} \begin{pmatrix}
y(0) \\ z(0)
\end{pmatrix},
\end{equation}
with non-trivial $(y(0),z(0))^t$ if and only if $\lambda$ is a root of the {\it characteristic equation}
\begin{equation}\label{sect6-eq9}
m_{11}(\lambda)m_{22}(\lambda)-m_{12}(\lambda)m_{21}(\lambda)=0,
\end{equation}
where, with $e_\lambda$ denoting the function defined by 
\begin{equation}\label{sect6-eq10}
e_\lambda(\theta):=e^{\lambda\theta},
\end{equation}
we have
\begin{equation}\label{sect6-eq11}
\begin{aligned}
m_{11}(\lambda)= & 1-D_1\, F_1(b_*,S_*)\,e_\lambda, \\
m_{12}(\lambda)= & - D_2\, F_1(b_*,S_*)\, e_\lambda \\
m_{21}(\lambda)= & -D_1\, F_2(b_*,S_*)\, e_\lambda \\
m_{22}(\lambda)= & \lambda-D_2\, F_2(b_*,S_*)\, e_\lambda.
\end{aligned}
\end{equation}

According to Section 5 of \cite{DG}, the steady state $(b_*,S_*)$ is asymptotically stable if all roots of \eqref{sect6-eq9} have negative real part, whereas it is unstable if there exists at least one root with positive real part. 

Combining \eqref{sect6-eq7} and \eqref{sect6-eq3} we deduce
\begin{equation}\label{sect6-eq12}
\begin{aligned}
D_1\, F_i(b_*,S_*)\,e_\lambda= & \int_0^\infty\delta\left(\overline{\xi}(a),S_*\right)\exp\left\{-\int_0^a\left(\lambda+\mu\left(\overline{\xi}(\tau),S_*\right)\right)\,\ud \tau\right\} \\
 = & \int_{x_b}^\infty\frac{\delta(x,S_*)}{g(x,S_*)}\exp\left\{-\int_{x_b}^x\frac{\lambda+\mu(y,S_*)}{g(y,S_*)}\,\ud y\right\}\,\ud x,
\end{aligned}
\end{equation}
where $\delta=\beta$ for $i=1$ and $\delta=-\gamma$ for $i=2$, and with $\overline{\xi}$ defined by \eqref{sect5-eq1}.

As elucidated by \eqref{sect6-eq3}, \eqref{sect6-eq2} and \eqref{X-s-def-2}, the dependence of $F_i$ on the $S$ variable involves the solution of the ODE describing how the size of an individual changes under the environmental condition described by $S$.  If we put
\begin{equation}\label{sect6-eq13}
X_{S_*+\varepsilon z_t}(0,-a,x_b)=\overline{\xi}(a)+\varepsilon\eta(a)+o(\varepsilon),
\end{equation}
then \eqref{X-s-def-2} implies that $\eta$ is a solution of
\begin{equation}\label{sect6-eq14}
\begin{aligned}
\dot{\eta}(\tau)= & D_1\, g\left(\overline{\xi}(\tau),S_*\right)\,\eta(\tau)+D_2\, g\left(\overline{\xi}(\tau),S_*\right)\, z_t(-a+\tau) \\
\eta(0)= & 0.
\end{aligned}
\end{equation}
It follows that 
\begin{equation}\label{sect6-eq15}
\eta(a)=\int_0^aK(a,\sigma)z_t(-a+\sigma)\,\ud \sigma,
\end{equation}
where
\begin{equation}\label{sect6-eq16}
K(a,\sigma):=D_2\, g\left(\overline{\xi}(\sigma),S_*\right)\exp\left\{\int_\sigma^aD_1\, g\left(\overline{\xi}(\theta),S_*\right)\ud\theta\right\}.
\end{equation}
Starting from \eqref{F-s-def-2} we find by straightforward Taylor expansion that 
\begin{align*}
 \mathcal{F}_{S_*+\varepsilon z_t}(0,-a,x_b) & = \mathcal{F}_{S_*}(0,-a,x_b) \\ 
& - \varepsilon\,\mathcal{F}_{S_*}(0,-a,x_b)\left(\int_0^a D_1\,\mu\left(\overline{\xi}(\tau),S_*\right)\eta(\tau)\,\ud\tau+\int_0^a D_2\,\mu\left(\overline{\xi}(\tau),S_*\right)z_t(-a+\tau)\,\ud \tau\right)+o(\varepsilon).
\end{align*}
For notational convenience we define
\begin{equation}\label{sect6-eq17}
\overline{\mathcal{F}}(a):= \exp\left\{-\int_0^a \mu\left(\overline{\xi}(\tau),S_*\right)\,\ud\tau\right\}=  \exp\left\{-\int_{x_b}^{\overline{\xi}(a)}\frac{\mu(x,S_*)}{g(x,S_*)}\,\ud x\right\}.
\end{equation}
Using the notation \eqref{sect6-eq17} we have
\begin{equation}\label{sect6-eq18}
\begin{aligned}
 D_2\, F_1(b_*,S_*)\, e_\lambda   = &   b_*\int_0^\infty D_2\, \beta\left(\overline{\xi}(a),S_*\right)\overline{\mathcal{F}}(a)\,\ud a + b_*\int_0^\infty D_1\, \beta\left(\overline{\xi}(a),S_*\right)\overline{\mathcal{F}}(a)\int_0^a K(a,\sigma) e^{\lambda(\sigma-a)}\,\ud\sigma\,\ud a \\
 - &  b_*\int_0^\infty  \beta\left(\overline{\xi}(a),S_*\right)\overline{\mathcal{F}}(a) \\
 &  \quad\times \left(\int_0^a D_1\, \mu\left(\overline{\xi}(\tau),S_*\right) \int_0^\tau K(\tau,\sigma)e^{\lambda(\sigma-a)}\,\ud\sigma\,\ud\tau+\int_0^a D_2\, \mu\left(\overline{\xi}(\tau),S_*\right) e^{\lambda(\tau-a)}\,\ud\tau\right)\,\ud a.
\end{aligned}
\end{equation}
The corresponding expression for $D_2\, F_2(b_*,S_*)\, e_\lambda$ is obtained from \eqref{sect6-eq18} by replacing $\beta$ by $\gamma$, multiplying the right hand side of \eqref{sect6-eq18} by $-1$ and adding $f'(S_*)$.

To illustrate the usefulness of the characteristic equation, we present an interesting instability result that, as far as we know, is new.

\begin{theorem}
A positive steady state $(b_*,S_*)$ of \eqref{DE-1} is unstable if $\mathcal{R}'(S_*)<0$ holds. 
\end{theorem}
{\bf Proof.} Let us denote the left hand side of \eqref{sect6-eq9} by $M(\lambda)$. We shall show that $M(0)<0$ and $M(\lambda)\to +\infty$ as $\lambda\to +\infty$. This then combined with the Intermediate Value Theorem implies that for some $\lambda>0$ we have $M(\lambda)=0$; and therefore the steady state is unstable.

First note that we have $m_{11}(0)=0$, since $D_1\,F_1(b_*,S_*)\, e_0=\mathcal{R}(S_*)=1$. Next observe that it follows from the definition of $m_{21}(\lambda)$ that $m_{21}(0)>0$. Since $F_1(b,S e_0) = b\,\mathcal{R}(S)$ for $S\geq 0$, differentiating with respect to $S$ (applying the chain rule in the left hand side) it follows $D_2F_1(b,S e_0)e_0 =b\,\mathcal{R}'(S)$, so that $m_{12}(0) = - b^*\,\mathcal{R}'(S_*)$ and hence $m_{12}(0) > 0$ if $\mathcal{R}'(S_*)<0$. Hence under this assumption we have $M(0) =-m_{12}(0) m_{21}(0) < 0$.

For $\lambda\to +\infty$ we have $D_1F_j(b_*,S_*)e_\lambda\to 0$ for $j=1,2$. Therefore,  $m_{21}(\lambda)\to 0$, $m_{11}(\lambda)\to 1$ and $m_{12}(\lambda)$ tends to a constant, while $m_{22}(\lambda)\sim \lambda$, 
and in particular this implies that $M(\lambda)\to \infty$, as $\lambda\to +\infty$; and the proof is complete. \eofproof

\begin{theorem}\label{Theorem 5.5}
Let $(n_*,S_*)$ be a steady state of the nonlinear semigroup $\mathcal{T}_{PDE}(t)$ defined in Theorem \ref{Theorem 2.4},  i.e., let $S_*$ be a positive root of (5.3) and let $n_*$ be defined by (5.9)-(5.10).
The assertions
\begin{itemize}
\item[(AS)] $(n_*,S_*)$ is locally asymptotically stable if the roots of the characteristic equation as specified in (5.21)-(5.30) are all in the left-half of the complex plane, at a uniform distance from the imaginary axis;
\item[(U)] $(n_*,S_*)$ is unstable if at least one root of this characteristic equation lies in the open right-half of the complex plane;
\end{itemize}
are true if the model ingredients $f, g, \mu, \beta$ and $\gamma$ are, such that the following hypotheses hold:
\begin{itemize}
\item[(i)] concerning the behaviour at ‘infinity’: H$_{g_\infty}$ ; H$_s$ ; $3\mu_0<\hat{\mu}$
(here $\mu_0$ is the weight that features in the definition of the state space and $\hat{\mu}$ is the asymptotic death rate);
\item[(ii)] concerning smoothness: H4$_f$ , H1$_h$ and H5$_h$ for $h=g, \mu, \beta, \gamma$;
\item[(iii)] strictly positive growth rate: H3$_g$;
\item[(iv)] boundedness: H2$_h$ for $h=g, \mu, \beta, \gamma$.
\end{itemize}
\end{theorem}
In essence, the proof is indirect, i.e., based on combining the corresponding result for the semigroup $\mathcal{T}_{DE}(t)$ with the continuity of the map $\mathcal{L}$ and its pseudo-inverse, cf. Theorem \ref{stab-eq-th}. It seems very likely that several of our assumptions can be relaxed. The problem of providing a direct proof is widely open.

Substituting \eqref{sect6-eq12}, \eqref{sect6-eq18} and its analogue for $D_2\, F_2$ into \eqref{sect6-eq11} we obtain a characteristic equation of the form \eqref{sect6-eq9} that is explicit in the ingredients of the model. In \cite{DGMNR} a more general variant of this characteristic equation was analysed in order to derive biological insight by unravelling the relationship between mechanisms at the individual level and phenomena (in particular oscillations) at the population level.

The continuous differentiability of $F_1, F_2$ is a sufficient condition for the differentiability of the nonlinear semigroup operators with respect to the initial state. It is not a necessary condition, see \cite{DK}. The more general variant of \eqref{sect6-eq9} derived in \cite{DGMNR} pertains to a model in which individual behaviour may change abruptly at the 
transition from juvenile to adult at a given size $\bar{x}$. For such a model, the maps $F_1,F_2$ are indeed not continuously differentiable exactly because of the state-dependent (i.e. food history dependent) delay between being born and becoming an adult (i.e. starting to reproduce). It is an open problem to prove the principle of linearised stability for this class of model.

\section{Concluding remarks}

While for age-structured population models there exists extensive literature (e.g., \cite{Webb1985,MagalRuan,DG}) justifying the Principle of Linearised Stability for steady states, there is as yet no such justification for general size-structured models. The reason is that such models are quasi-linear, in the sense that the nonlinearity affects the highest derivative. Concerning special models, we are aware of  \cite{Grabosch1990,Grabosch1991}, also see \cite{Grabosch1989}, in which a separable growth rate $g(x,S)=g_1(x)g_2(S)$ is assumed, allowing an implicit time transformation that, in a sense, eliminates the nonlinearity.

Here we concentrated on the so-called Daphnia model, in which all newborns are assumed to have a fixed given size $x_b$ and the nonlinearity is due to competition for food. A consequence of the fixed birth size is that the population dynamics is ‘driven’ by a scalar renewal equation. By working with the history of the birth rate, rather than the current size distribution, we obtain a delay equation formulation of the problem. The advantage is that the  corresponding dynamics is based on translation of information with fixed (rather than variable) speed, allowing rigorous linearisation.

In order to assess (in)stability in terms of size distributions, we have studied the precise relationship between the ‘current size’ and the ‘age + history of food’ ways of bookkeeping. Thus we were able to transfer stability information from one framework to the other. As far as our literature search revealed, this is a new approach leading rather indirectly to new PDE results.

What next? For general size-structured models one can work with a renewal equation for a function taking values in an infinite dimensional space. This leads, as far as we know, into unexplored territory (but see \cite{FrancoGyllenbergDiekmann} for strong results under very restrictive assumptions).

\appendix

\section{Proof of Lemma \ref{lemma2.4}}

\noindent By hypothesis there exist constants $B, \Gamma, G_1, G_2$, $M_1, M_2, B_1, B_2, \Gamma_1, \Gamma_2$ and $F_1$ such that
\begin{equation}\label{beta-gamma-bounds}
\begin{aligned}
|\beta(x,S)|\leq & \,B,\\
|\gamma(x,S)|\leq & \,\Gamma
,
\end{aligned}
\end{equation}
and
\begin{equation}\label{ingredient-lipschitzCond}
\begin{aligned}
|g(x_1,S_1)-g(x_2,S_2)|\leq & G_1|x_1-x_2|+G_2|S_1-S_2|,\\
|\mu(x_1,S_1)-\mu(x_2,S_2)|\leq & M_1|x_1-x_2|+M_2|S_1-S_2|,\\
|\beta(x_1,S_1)-\beta(x_2,S_2)|\leq & B_1|x_1-x_2|+B_2|S_1-S_2|,\\
|\gamma(x_1,S_1)-\gamma(x_2,S_2)|\leq & \Gamma_1|x_1-x_2|+\Gamma_2|S_1-S_2|,\\
|f(S_1)-f(S_2)|\leq & F_1|S_1-S_2|,
\end{aligned}
\end{equation}
which imply, using Gr\"{o}nwall's inequality (see  definitions (\ref{X-s-def-2}), (\ref{F-s-def-2}),  (\ref{RE-3}) and (\ref{gamma-2})),
\begin{equation}\label{ingredient-bounds}
\begin{aligned}
|X_{S_1}(t,s,\xi)-X_{S_2}(t,s,\xi)|\leq & G_2 e^{G_1(t-s)}\sup_{\tau\in[s,t]}|S_1(\tau)-S_2(\tau)|,\\
|\mathcal{F}_{S_1}(t,s,\xi)-\mathcal{F}_{S_2}(t,s,\xi)|\leq & (t-s)\left(M_1G_2 e^{G_1(t-s)}+M_2\right)\sup_{\tau\in[s,t]}|S_1(\tau)-S_2(\tau)|,\\
|\beta_{S_1}(t,s,\xi)-\beta_{S_2}(t,s,\xi)|\leq & H_1(t-s)\sup_{\tau\in[s,t]}|S_1(\tau)-S_2(\tau)|,\\
|\gamma_{S_1}(t,s,\xi)-\gamma_{S_2}(t,s,\xi)|\leq & \tilde{H}_1(t-s)\sup_{\tau\in[s,t]}|S_1(\tau)-S_2(\tau)|,\\
\end{aligned}
\end{equation}
where 
\begin{equation}
\begin{aligned}
H_1(t):= & B_1G_2e^{G_1t}+B_2+Bt\left(M_1 G_2 e^{G_1t}+M_2\right),\\
\tilde{H}_1(t):= & \Gamma_1G_2e^{G_1t}+\Gamma_2+\Gamma t\left(M_1 G_2 e^{G_1 t}+M_2\right).
\end{aligned}
\end{equation}
From the third bound in (\ref{ingredient-bounds}) and definition (\ref{f-def}) it follows
\begin{equation}\label{h-bound}
|h_{S_1}(t)-h_{S_2}(t)|\leq H_1(t)\|n_0\|_{L^1}\sup_{\tau\in[0,t]}|S_1(\tau)-S_2(\tau)|,
\end{equation}
and, analogously, from the fourth bound in (\ref{ingredient-bounds}) and definition (\ref{ks-def}) it follows
\begin{equation}\label{k-diff-bound}
|k_{S_1}(t)-k_{S_2}(t)|\leq \tilde{H}_1(t)\|n_0\|_{L^1}\sup_{\tau\in[0,t]}|S_1(\tau)-S_2(\tau)|.
\end{equation}
From \eqref{RE} and \eqref{f-def}  it follows
\begin{equation}
b_{S}(t)\leq B \int_0^t b_S(t-a)\ud a + B\|n_0\|_{L^1},
\end{equation}
so that, using Gr\"{o}nwall's inequality,
\begin{equation}\label{birthrate-bound}
b_S(t)\leq B \|n_0\|_{L^1} e^{Bt}.
\end{equation}
From \eqref{ingredient-bounds}, \eqref{h-bound} and \eqref{birthrate-bound} we have 
\begin{equation}\label{birth-diff-bound}
\begin{aligned}
|b_{S_1}(t)-b_{S_2}(t)|\leq & \int_0^t |\beta_{S_1}(t,\tau,x_b)-\beta_{S_2}(t,\tau,x_b)|b_{S_1}(\tau)\ud\tau+\\
&\quad+\int_0^t \beta_{S_2}(t,\tau,x_b)|b_{S_1}(\tau)-b_{S_2}(\tau)|\ud\tau+|h_{S_1}(t)-h_{S_2}(t)|\\
\leq & H_2(t)\|n_0\|_{L^1}\sup_{\tau\in[0,t]}|S_1(\tau)-S_2(\tau)|+B\int_0^t |b_{S_1}(\tau)-b_{S_2}(\tau)| \ud\tau
\end{aligned}
\end{equation}
with
\begin{equation}
H_2(t)=B \int_0^tH_1(t-\tau)e^{B\tau}\ud\tau+H_1(t).
\end{equation}
Using that $H_2(t) e^{Bt}||n_0||_{L^1}\displaystyle\sup_{\tau\in [0,t]}|S_1(\tau)-S_2(\tau)|$  is an increasing function of $t$ we obtain from Gr\"{o}nwall's inequality
\begin{equation}\label{birth-diff-bound_2}
|b_{S_1}(t)-b_{S_2}(t)|\leq H_2(t)e^{Bt}\|n_0\|_{L^1}\sup_{\tau\in[0,t]}|S_1(\tau)-S_2(\tau)|.
\end{equation}
Now, from \eqref{volterraeq} we can bound using \eqref{birthrate-bound}, \eqref{ingredient-lipschitzCond}, \eqref{ingredient-bounds}, \eqref{birth-diff-bound_2} and \eqref{k-diff-bound},
\begin{equation}
\begin{aligned}
|\mathcal{V}(S_1)(t)-\mathcal{V}(S_2)(t)|\leq & \int_0^t |f(S_1(\tau))-f(S_2(\tau))|d\tau +\\
& + B \int_0^t\int_0^\tau|\gamma_{S_1}(\tau,s,x_b)-\gamma_{S_2}(\tau,s,x_b)| e^{Bs}\ud s\ud\tau \|n_0\|_{L^1}+\\
& +\Gamma\int_0^t\int_0^\tau|b_{S_1}(s)-b_{S_2}(s)|\ud s+\int_0^t |k_{S_1}(\tau)-k_{S_2}(\tau)|\ud\tau\leq\\
\leq & \int_0^t L(\tau)\sup_{s\in[0,\tau]} |S_1(s)-S_2(s)|d\tau
\end{aligned}
\end{equation}
with
\begin{equation}
L(\tau) := F_1+B\int_0^\tau\tilde{H}_1(\tau-s)e^{Bs} \ud s\|n_0\|_{L^1}+\Gamma\int_0^\tau H_2(s)e^{Bs}\ud s \|n_0\|_{L^1} + \tilde{H}_1(\tau)\|n_0\|_{L^1}.
\end{equation}
The function $L$ is increasing, so that for $t\in[0,\hat{t}]$ and any $k>0$ we have 
\begin{equation}
\begin{aligned}
e^{-kt}|\mathcal{V}(S_1)(t)-\mathcal{V}(S_2)(t)|\leq & L(\hat{t})\int_0^t 
e^{-k(t-\tau)}e^{-k\tau}\sup_{s\in[0,\tau]} |S_1(s)-S_2(s)|\ud\tau\\
\leq & L(\hat{t})\int_0^t e^{-k(t-\tau)}\left(\sup_{s\in[0,\hat{t}]} e^{-k s} |S_1(s)-S_2(s)|\right)\ud\tau\\
\leq & \frac{L(\hat{t})}{k}\|S_1-S_2\|_{C};
\end{aligned}
\end{equation}
where we have defined, on the space $C\left([0,\hat{t}],\mathbb{R}\right)$, the norm 
\begin{equation}\label{norm-def}
||S||_C:=\displaystyle\sup_{s\in [0,\hat{t}]}e^{-ks}|S(s)|. 
\end{equation}
We have
\begin{equation}\label{contraction}
\|\mathcal{V}(S_1)-\mathcal{V}(S_2)\|_C=\sup_{t\in[0,\hat{t}]} e^{-kt}|\mathcal{V}(S_1)(t)-\mathcal{V}(S_2)(t)| \leq \frac{L(\hat{t})}{k}\|S_1-S_2\|_C.
\end{equation}
Thus for any fixed $\hat{t}$ we choose $k$ such that $L(\hat{t})<k$ holds, and therefore $\mathcal{V}$ is a contraction on the Banach-space $C([0,\hat{t}],\mathbb{R})$ endowed with the norm $\|\cdot\|_C$ defined in \eqref{norm-def}.

\section{Continuity of $\mathcal{L}$ and $\mathcal{L}_{ps}^{-1}$}\label{appendix-continuityL}

\noindent In this appendix the letters $\sigma$ and $\alpha$ are used as mnemonic labels to the words ``size'' and ``age'' respectively. Let $\tilde{k}>0$ be such that  $-\tilde{k} < \partial_1 g < \tilde{k}$ and define $k_\alpha=\tilde{k}+1/a_r$ and $k_\sigma=g_{min}^{-1}\tilde{k}+1/x_r$ where $a_r>0$ and $x_r>0$ are a referential age and a referential size, respectively that can be chosen freely (notice that $k_\alpha^{-1}$ has time units and $k_\sigma^{-1}$ has size units). For each $h\in C^1([x_b,\infty))$ define the size weighted norm 
$$
|h|_\sigma=\sup_{x\in[x_b,\infty)} |h(x)e^{-\frac{x}{x_r}}|+\sup_{x\in[x_b,\infty)}, 
|h'(x)e^{-k_\sigma x}|
$$
and the subset of $C^1([x_b,\infty))$ given by
$$
W_\sigma=\left\{h\in C^1([x_b,\infty)), \quad\text{such that}\quad h([x_b,\infty))= (-\infty,0]\text{, } h'(x)< 0 \text{ and }|h|_\sigma <\infty \right\}.
$$
Similarly, for $h\in C^1((-\infty,0])$ define the age weighted norm
$$
|h|_\alpha=\sup_{a\in[0,\infty)} |h(-a)e^{-\frac{a}{a_r}}|+\sup_{a\in[0,\infty)} |h'(-a)e^{-k_\alpha a }|,
$$
and the subset of $C^1((-\infty,0])$ given by
$$
W_\alpha=\left\{h\in C^1((-\infty,0]) \quad\text{such that}\quad h((-\infty,0])= [x_b,\infty)\text{, } h'(-a) < 0 \text{ and }|h|_\alpha <\infty \right\}.
$$
Let 
$$W^{g_{\infty}}_\sigma:=\left\{h\in W_\sigma | h(x) \leq \frac{c_2 - x}{g_\infty}\right\},$$ 
with the subspace topology ($W^{g_{\infty}}_\sigma \subset W_\sigma$) and 
$$W^{g_{\infty}}_\alpha:=\left\{h\in W_\alpha | c_1+g_\infty a \leq h(-a)\right\},$$ 
with the subspace topology ($W^{g_{\infty}}_\alpha \subset W_\alpha$), where $c_1$ and $c_2$ are the constants given in Lemma \ref{lemmaBoundsXandT}, i.e.
\[
c_1=\bar{x}-\frac{\bar{x}-x_b}{g_\text{min}}g_\infty,\qquad\text{and}\qquad c_2=\bar{x}.
\]

\begin{lemma}\label{cL_XandTauCont}
The mappings
$$
\begin{array}{cccl}
T:&\X_2 & \longrightarrow & W^{g_{\infty}}_\sigma\\
&\psi & \longmapsto & T(\psi):=T_\psi(x_b,\cdot,0),\\
\end{array}
$$
and
$$
\begin{array}{cccl}
X:&\X_2 & \longrightarrow & W^{g_{\infty}}_\alpha\\
&\psi & \longmapsto & X(\psi):=X_\psi(0,\cdot,x_b),\\
\end{array}
$$
are well defined and continuous.
\end{lemma}

{\bf Proof.} Recall that $X_\psi(t,s,x_b)$ is the solution of
\begin{equation}
\left\{
\begin{aligned}
\partial_t X_\psi(t,s,x_b)= &g(X_\psi(t,s,x_b),\psi(t))\\
X_\psi(s,s,x_b)=& x_b
\end{aligned}
\right. .
\end{equation}
Since $0<g_{min}\leq g \leq g_{max}$, it follows that $X(\psi)(-a)\in[x_b+g_{min}a,x_b+g_{max}a]$, so that $X(\psi)(0)=x_b$ and $\displaystyle\lim_{a\rightarrow \infty} X(\psi)(-a)=\infty$, which implies $X(\psi)((-\infty,0])=[x_b,\infty)$. An expression for $X(\psi)'(-a)$ is obtained through the variational equation of the above initial value problem. Specifically we consider the initial value problem obtained for $\partial_s X_\psi(t,s,x_b)$, that is
\[
\left\{
\begin{aligned}
\partial_t \left(\partial_2 X_\psi(t,s,x_b)\right)= &\partial_1 g(X_\psi(t,s,x_b),\psi(t))\partial_2 X_\psi(t,s,x_b)\\
\partial_t X_\psi(s,s,x_b)+\partial_2 X_\psi(s,s,x_b)= & 0
\end{aligned}
\right. ,
\]
whose solution is
\[
\partial_2 X_\psi(t,s,x_b)=-g(x_b,\psi(s))e^{\int_{s}^t \partial_1 g(X_\psi(\tau,s,x_b),\psi(\tau))d\tau}.
\]
Then, since $X(\psi)'(-a)=\partial_2 X_\psi(0,-a,x_b)$, we have
$$
X(\psi)'(-a)=-g(x_b,\psi(-a))e^{\int_{-a}^0 \partial_1 g(X_\psi(\tau,-a,x_b),\psi(\tau))d\tau}.
$$
Using that $-\tilde{k} < \partial_1 g < \tilde{k}$ for some $\tilde{k}>0$, it follows $$
-g_{max} e^{a\tilde{k}}< X(\psi)'(-a)< -g_{min} e^{-a \tilde{k}}.
$$
Then, since $k_\alpha=\tilde{k}+1/a_r$, the bounds on $X(\psi)$ and $X(\psi)'$ imply that $X(\psi)\in W_\alpha$, and then the lower bound in (\ref{X-upper-bound}) implies $ X(\psi)\in W^{g_{\infty}}_\alpha$.

Since $T(\psi)=T_\psi(x_b,\cdot,0)$ is the inverse function of $X(\psi)=X_\psi(0,\cdot,x_b)$, it follows that $T(\psi)(x)\in [(x_b-x)g_{min}^{-1},(x_b-x)g_{max}^{-1}]$, so that $T(\psi)(x_b)=0$ and $\displaystyle\lim_{x\rightarrow \infty} T(\psi)(x)=-\infty$, which implies $T(\psi)([x_b,\infty))=(-\infty,0]$. Moreover, taking into account the bound on the derivative of $X(\psi)$, one obtains
$$
-g_{min}^{-1} e^{-T(\psi)(x) \tilde{k}}<T(\psi)'(x)<-g_{max}^{-1} e^{T(\psi)(x) \tilde{k}},
$$
so that, using the bounds on $T(\psi)$,
$$
-g_{min}^{-1} e^{(x-x_b)g_{min}^{-1} \tilde{k}}<T(\psi)'(x)<-g_{max}^{-1} e^{(x_b-x)g_{min}^{-1} \tilde{k}}.
$$
Then, since $k_\sigma=g_{min}^{-1} \tilde{k}+1/x_r$, the bounds on $T(\psi)$ and $T(\psi)'$ imply that $T(\psi)\in W_\sigma$, and then the upper bound in (\ref{tau-upper-bound}) implies $ T(\psi)\in W^{g_{\infty}}_\sigma$.

To prove the continuity of $X$ with respect $\psi$, let a sequence $\{\psi_n\}_{n\in\mathbb{N}}$ converge to $\psi$ and let us show that $|X(\psi_n)-X(\psi)|_\alpha<\varepsilon$ for all $n$ large enough. Indeed, take $\bar{a}>0$, such that
$$
\sup_{a\in(\bar{a},\infty)} g_{max}\,a e^{- \frac{a}{a_r}} <\frac{\varepsilon}{4},\quad\text{ and }\quad 
\sup_{a\in(\bar{a},\infty)} g_{max}\, e^{-\frac{a}{a_r}} <\frac{\varepsilon}{4}.
$$
Next, take $n$ large enough so that
$$
\sup_{a\in[0,\bar{a}]}|X(\psi_n)(-a)-X(\psi)(-a)| < \frac{\varepsilon}{4},\quad\text{ and }\quad
\sup_{a\in[0,\bar{a}]}|X(\psi_n)'(-a)-X(\psi)'(-a)| < \frac{\varepsilon}{4},
$$
which is possible because the sequence $\psi_n$ converges uniformly to $\psi$ on $[-\bar{a},0]$, so that the sequences $X(\psi_n)$ and $X(\psi_n)'$ converge also uniformly on $[-\bar{a},0]$. Then, it follows
$$
\begin{aligned}
&\sup_{a\in[0,\infty)}|X(\psi_n)(-a)-X(\psi)(-a)|e^{- \frac{a}{a_r}}\\
\leq & \sup_{a\in[0,\bar{a}]}|X(\psi_n)(-a)-X(\psi)(-a)|e^{- \frac{a}{a_r}} + \sup_{a\in(\bar{a},\infty)}g_{max}a e^{- \frac{a}{a_r}}<\frac{\varepsilon}{2},
\end{aligned}
$$
and
$$
\begin{aligned}
&\sup_{a\in[0,\infty)}|X(\psi_n)'(-a)-X(\psi)'(-a)|e^{-k_\alpha a} =\sup_{a\in[0,\infty)}|X(\psi_n)'(-a)-X(\psi)'(-a)|e^{-\left(\tilde{k}+\frac{1}{a_r}\right) a} \\
\leq & \sup_{a\in[0,\bar{a}]}|X(\psi_n)'(-a)-X(\psi)'(-a)|e^{-\left(\tilde{k}+\frac{1}{a_r}\right) a} + \sup_{a\in(\bar{a},\infty)}g_{max} e^{a\tilde{k}} e^{-\left(\tilde{k}+\frac{1}{a_r}\right)a}<\frac{\varepsilon}{2},
\end{aligned}
$$
if $n$ is large enough, so that $|X(\psi_n)-X(\psi)|_\alpha<\varepsilon$ for these values of $n$. 
An analogous argument is made to prove the continuity of $T$ with respect $\psi$.
\eofproof

\begin{lemma}\label{cL_changeVar}
The mappings
$$
\begin{array}{cccl}
\zeta_\sigma:&L^1_{-(\hat{\mu}-\mu_0)}((-\infty,0],\mathbb{R}_+)\times W^{g_{\infty}}_\sigma & \longrightarrow & L^1_{\kappa_0}([x_b,\infty),\mathbb{R}_+)\\
&(\eta,h) & \longmapsto & -\eta(h(\cdot))h'(\cdot)\\
\end{array}
$$
and
$$
\begin{array}{cccl}
\zeta_\alpha:&L^1_{\kappa_0}([x_b,\infty),\mathbb{R}_+)\times W^{g_{\infty}}_\alpha & \longrightarrow & L^1_{-(\hat{\mu}-\mu_0)}((-\infty,0],\mathbb{R}_+)\\
&(\eta,h) & \longmapsto & -\eta(h(\cdot))h'(\cdot)\\
\end{array}
$$
are well defined and continuous.
\end{lemma}
{\bf Proof.} 
To see that $\zeta_\sigma$ is well defined notice, using $h'(x)<0$ and $\eta(y)\ge 0$, that
$$
\int_{x_b}^{\infty}|\eta(h(x))h'(x)|e^{\kappa_0 x}dx =-\int_{x_b}^{\infty}\eta(h(x))h'(x)e^{\kappa_0 x}dx = -\int_{h(x_b)}^{h(\infty)}\eta(y)e^{\kappa_0 h^{-1}(y)}dy = \int_{-\infty}^{0}\eta(y)e^{\kappa_0 h^{-1}(y)}dy,
$$
and then, since $h(x)\leq (\bar{x}-x)g^{-1}_{\infty}$ implies that $h^{-1}(y)\leq \bar{x}-g_{\infty}y$ and $\kappa_0 g_{\infty}=\hat{\mu}-\mu_0$ by definition,
$$
\int_{-\infty}^{0}\eta(y)e^{\kappa_0 h^{-1}(y)}dy \leq e^{\kappa_0 \bar{x}}\int_{-\infty}^{0}\eta(y) e^{-(\hat{\mu}-\mu_0)y}dy <\infty.
$$

To prove that $\zeta_\sigma$ is continuous we take a sequence $\{(\eta_n,h_n)\}_{n\in\mathbb{N}}$ converging to $(\eta,h)$ and we show that $\|\zeta_\sigma(\eta_n,h_n)-\zeta_\sigma(\eta,h)\|<\varepsilon$ for $n$ large enough.

To do so, first let $\bar{x}\in[x_b,\infty)$ be such that
$$
e^{\kappa_0 x_b}\int_{-\infty}^{h(\bar{x})+1}\eta(y) e^{-(\hat{\mu}-\mu_0)}dy < \varepsilon/4,
$$
which is possible because $\displaystyle\lim_{x\rightarrow\infty} h(x)=-\infty$ and $\eta\in L^1_{-(\hat{\mu}-\mu_0)}$.

Then consider
$$
\|\zeta_\sigma(\eta_n,h_n)-\zeta_\sigma(\eta,h)\|=\int_{x_b}^\infty |\eta_n(h_n(x))h_n'(x)-\eta(h(x))h'(x)|e^{\kappa_0 x} dx \leq I_{1} + I_{21} 
+ + I_{22} + I_{23}
$$
with
$$
\begin{aligned}
I_{1}& = \int_{x_b}^{\infty} |\eta_n(h_n(x))h_n'(x)-\eta(h_n(x))h_n'(x)|e^{\kappa_0 x}dx \\
& = -\int_{x_b}^{\infty} |\eta_n(h_n(x))-\eta(h_n(x))|h_n'(x) e^{\kappa_0 x} dx = \int_{-\infty}^{h_n(x_b)} |\eta_n(y)-\eta(y)| e^{\kappa_0 h_n^{-1}(y)}dy \\
& \leq \int_{-\infty}^{0} |\eta_n(y)-\eta(y)| e^{\kappa_0 (\bar{x}-g_{\infty} 
y)}dy \leq e^{\kappa_0 \bar{x}} \int_{-\infty}^{0} |\eta_n(y)-\eta(y)| e^{ -(\hat{\mu}-\mu_0) y}dy \leq e^{\kappa_0 \bar{x}} \|\eta_n-\eta\|_{L^1_{-(\hat{\mu}-\mu_0)}}\\
I_{21}& = \int_{x_b}^{\bar{x}} |\eta(h_n(x))h_n'(x)-\eta(h(x))h'(x)| e^{\kappa_0 x}dx, \\
I_{22}& = \int_{\bar{x}}^{\infty} |\eta(h_n(x))h_n'(x)| e^{\kappa_0 x} dx = -\int_{\bar{x}}^{\infty} \eta(h_n(x))h_n'(x)e^{\kappa_0 x}dx = \int_{-\infty}^{h_n(\bar{x})} \eta(y)e^{\kappa_0 h_n^{-1}(y)}dy \\
& \leq e^{\kappa_0 \bar{x}}\int_{-\infty}^{h_n(\bar{x})} \eta(y)e^{-(\hat{\mu}-\mu_0)y}dy\\
I_{23}& = \int_{\bar{x}}^{\infty} |\eta(h(x))h'(x)|e^{\kappa_0 x}dx = 
-\int_{\bar{x}}^{\infty} \eta(h(x))h'(x)e^{\kappa_0 x}dx = \int_{-\infty}^{h(\bar{x})} \eta(y)e^{\kappa_0 h^{-1}(y)}dy\\ & \leq e^{\kappa_0 \bar{x}}\int_{-\infty}^{h(\bar{x})} \eta(y)e^{-(\hat{\mu}-\mu_0)y}dy < \frac{\varepsilon}{4}.
\end{aligned}
$$
Clearly, for $n$ large enough $I_{1}<\varepsilon/4$ because $\eta_n\rightarrow \eta$. Since $h_n(\bar{x})\rightarrow h(\bar{x})$, it follows $ h_n(\bar{x}) <h(\bar{x}) + 1$ if $n$ is large enough, so that in these cases 
$I_{22}<\frac{\varepsilon}{4}$. Finally $I_{21}$ can also be bounded by $\varepsilon/4$ if $n$ is large enough because, as $n\rightarrow \infty$, $h_n\rightarrow h$ and $h_n'\rightarrow h'$ uniformly on the closed interval $[x_b,\bar{x}]$. This is analogous to the continuity of translation in $L^1$, which is a consequence of the density in $L^1$ of the space of continuous compactly supported functions. Indeed, one can choose $\theta$ continuous on a compact interval contained in $(-\infty, 0]$ 
such that
$$
e^{\kappa_0 \bar{x}}\int_{-\infty}^0 |\eta(y)-\theta(y)|e^{-(\hat{\mu}-\mu_0)y}dy < \frac{\varepsilon}{16}.
$$ 
Then, using the same changes of variables as before we have
$$
\begin{aligned}
I_{21}& = \int_{x_b}^{\bar{x}} |\eta(h_n(x))h_n'(x)-\eta(h(x))h'(x)| e^{\kappa_0 x}dx \leq \int_{x_b}^{\bar{x}} |\eta(h_n(x))h_n'(x)-\theta(h_n(x))h_n'(x)|e^{\kappa_0 x}dx   \\
& + \int_{x_b}^{\bar{x}} |\theta(h_n(x))h_n'(x)-\theta(h(x))h'(x)|e^{\kappa_0 x}dx + \int_{x_b}^{\bar{x}} |\theta(h(x))h'(x)-\eta(h(x))h'(x)|e^{\kappa_0 x}dx\\
& \leq 2 e^{\kappa_0 \bar{x}} \int_{-\infty}^{0} |\eta(y)-\theta(y)|e^{-(\hat{\mu}-\mu_0)y} dy + \int_{x_b}^{\bar{x}} |\theta(h_n(x))h_n'(x)-\theta(h(x))h'(x)|e^{\kappa_0 x}dx  < \frac{\varepsilon}{4},
\end{aligned}
$$
when $n$ is large enough due to the uniform convergence of $h_n$ and $h_n'$ and the uniform continuity of $\theta$ which imply that the second term is smaller that $\varepsilon/8$ for these large $n$.

Analogous arguments can be followed to show that $\zeta_\alpha$ is also well defined and continuous.  

\eofproof

\begin{lemma}\label{cL_prod}
Let $\mathcal{X}_{\mu_0}=\mathcal{X}_{1,\mu_0}\times\mathcal{X}_{2,\mu_0}$ be as in Definition \ref{def-phase-space-DE}. Then the mappings
$$
\begin{array}{cccl}
p_\sigma:&\X_{1,\mu_0}\times\X_{2,\mu_0} & \longrightarrow & L^1_{-(\hat{\mu}-\mu_0)}((-\infty,0],\mathbb{R}_+)\\
&(\phi,\psi) & \longmapsto & \phi(\cdot)\mathcal{F}_\psi(0,\cdot,x_b)\\
\end{array}
$$
and
$$
\begin{array}{cccl}
p_\alpha:&L^1_{-(\hat{\mu}-\mu_0)}((-\infty,0],\mathbb{R}_+)\times \X_{2,\mu_0} & \longrightarrow & \X_{1,\mu_0}\\
&(\eta,\psi) & \longmapsto & \frac{\eta(\cdot)}{\mathcal{F}_\psi(0,\cdot,x_b)}\\
\end{array}
$$
are well defined and continuous.
\end{lemma}
{\bf Proof.} Assumption (\ref{weight-ass1}) on $\mathcal{F}$ guarantees that both $p_\sigma$ and $p_\alpha$ are well defined. Indeed, $p_\sigma$ is well defined because 
$$
\int_{-\infty}^{0}\phi(\theta)\mathcal{F}_\psi(0,\theta,x_b) e^{-(\hat{\mu}-\mu_0)\theta}d\theta \leq C \int_{-\infty}^{0}\phi(\theta) e^{\hat{\mu}\theta} e^{-(\hat{\mu}-\mu_0)\theta}d\theta \leq C \int_{-\infty}^{0}\phi(\theta) e^{\mu_0 \theta}d\theta <\infty,
$$
and $p_\alpha$ is well defined because
$$
\int_{-\infty}^{0}\frac{\eta(\theta)}{\mathcal{F}_\psi(0,\theta,x_b)} e^{\mu_0 \theta}d\theta \leq \frac{1}{c} \int_{-\infty}^{0}\eta(\theta) e^{-\hat{\mu}\theta} e^{\mu_0\theta}d\theta \leq \frac{1}{c} \int_{-\infty}^{0}\eta(\theta) e^{-(\hat{\mu}-\mu_0) \theta}d\theta <\infty.
$$
To prove the continuity of $p_\sigma$ we take a sequence $\{(\phi_n,\psi_n)\}_{n\in\mathbb{N}}$ converging to $(\phi,\psi)$ and we show $\|p_\sigma(\phi_n,\psi_n)-p_\sigma(\phi,\psi)\|<\varepsilon$ for $n$ large enough. 
In the following we write $\mathcal{F}_\psi(\cdot)$ instead of $\mathcal{F}_\psi(0,\cdot,x_b)$ so that notation is simplified. Let $\bar{\theta}$ such that
$$
C \int_{-\infty}^{\bar{\theta}}  e^{\mu_0 \theta}|\phi(\theta)|d\theta <\frac{\varepsilon}{4}.
$$
Consider
$$
\|p_\sigma(\phi_n,\psi_n)-p_\sigma(\phi,\psi)\|=\int_{-\infty}^0 |\mathcal{F}_{\psi_n}(\theta)\phi_n(\theta)-\mathcal{F}_{\psi}(\theta)\phi(\theta)|e^{-(\hat{\mu}-\mu_0)\theta} d\theta \leq I_{11} + I_{12} + I_{21} + I_{22}
$$
with
$$
\begin{aligned}
&I_{11} = \int_{\bar{\theta}}^0 |\mathcal{F}_{\psi_n}(\theta)-\mathcal{F}_{\psi}(\theta)||\phi_n(\theta)|e^{-(\hat{\mu}-\mu_0)\theta}d\theta\leq 
\|\phi_n\|_{1,\mu_0}\sup_{\theta\in[\bar{\theta},0]} |\mathcal{F}_{\psi_n}(\theta)-\mathcal{F}_{\psi}(\theta)|e^{-\hat{\mu} \theta}, \\
&I_{12} = \int_{\bar{\theta}}^0 |\mathcal{F}_{\psi}(\theta)||\phi_n(\theta)-\phi(\theta)|e^{-(\hat{\mu}-\mu_0)\theta}d\theta \leq C \|\phi_n-\phi\|_{1,\mu_0}, \\
&I_{21} = \int_{-\infty}^{\bar{\theta}} |\mathcal{F}_{\psi}(\theta)||\phi_n(\theta)-\phi(\theta)|e^{-(\hat{\mu}-\mu_0)\theta}d\theta \leq C \|\phi_n-\phi\|_{1,\mu_0}, \\
&I_{22}= \int_{-\infty}^{\bar{\theta}} |\mathcal{F}_{\psi_n}(\theta)-\mathcal{F}_{\psi}(\theta)||\phi(\theta)|e^{-(\hat{\mu}-\mu_0)\theta} d\theta \leq \int_{-\infty}^{\bar{\theta}} Ce^{\mu_0 \theta}|\phi(\theta)|d\theta<\frac{\varepsilon}{4}.
\end{aligned}
$$
Then, for $n$ large enough $I_{12}$ and $I_{21}$ are smaller than $\varepsilon/4$ because $\phi_n\rightarrow \phi$ and $I_{11}$ is smaller than $\varepsilon/4$ because the sequence $\psi_n$ converges uniformly to $\psi$ 
on $[\bar{\theta},0]$ (so does the sequence $\mathcal{F}(\psi_n)$ to $\mathcal{F}(\psi)$) and $\|\phi_n\|_{1,\mu_0}$ stays bounded. So, for these values of $n$ one has $\|p_\sigma(\phi_n,\psi_n)-p_\sigma(\phi,\psi)\|<\varepsilon$.

Similar arguments can be followed to show that $p_\alpha$ is also continuous.
\eofproof

\begin{theorem} (Theorem \ref{contiuity-L}, Theorem \ref{pseudo-th}) Assume H$_s$, H2$_g$, H3$_g$ and H$_{g_\infty}$. Then the maps $\mathcal{L}$ and $\mathcal{L}^{-1}_{ps}$ are continuous.
\end{theorem}
{\bf Proof.} First consider $\mathcal{L}:\mathcal{X}_{\mu_0}\to L^1_{\kappa_0}\times\mathbb{R}$. The second component of $\mathcal{L}$ given by $\mathcal{L}_2(\phi,\psi)=\psi(0)$ is clearly continuous. The first component of this mapping can be written as the composition
$$
\begin{array}{cccccl}
\mathcal{L}_1:&\X_{1,\mu_0}\times\X_{2,\mu_0} & \longrightarrow & L^1_{-(\hat{\mu}-\mu_0)}(\mathbb{R}_-,\mathbb{R}_+)\times W^{g_{\infty}}_\sigma & \longrightarrow & L_{\kappa_0}^1([x_b,\infty),\mathbb{R}_+)\\
& (\phi,\psi) & \longmapsto & (p_\sigma(\phi,\psi),T(\psi)) & \longmapsto & \zeta_\sigma(p_\sigma(\phi,\psi),T(\psi)) \\
\end{array},
$$
where $\zeta_\sigma$ is the first mapping given in Lemma \ref{cL_changeVar}. The first function of this composition is continuous because of Lemmas \ref{cL_XandTauCont} and \ref{cL_prod}. The second function is continuous because of Lemma \ref{cL_changeVar}. Therefore, $\mathcal{L}$ is continuous.

Now consider $\mathcal{L}^{-1}_{ps}:L^1_{\kappa_0}\times\mathbb{R}\to \mathcal{X}_{\mu_0}$. Its second component is given by a bounded and linear operator, namely $\mathcal{L}^{-1}_{ps,2}(n,S)=S e_0$ where $e_0$ is the constant function equal to 1 (defined in $\mathbb{R}_-$), so that in particular $\mathcal{L}^{-1}_{ps,2}$ is continuous. The first component of $\mathcal{L}^{-1}_{ps}$  can be decomposed as
$$
\begin{array}{cccccl}
\mathcal{L}_{ps,1}^{-1}:&L^1_{\kappa_0}([x_b,\infty),\mathbb{R}_+)\times\mathbb{R}_+ & \longrightarrow & L^1_{-(\hat{\mu}-\mu_0)}(\mathbb{R}_-,\mathbb{R}_+)\times \X_{2,\mu_0} & \longrightarrow & \X_{1,\mu_0}\\
& (n,S) & \longmapsto & (\zeta_\alpha(n,X(S e_0)),S e_0) & \longmapsto & p_\alpha(\zeta_\alpha(n,X(S e_0)),S e_0) \\
\end{array},
$$
where $\zeta_\alpha$ is the second mapping given in Lemma \ref{cL_changeVar}. The second function of this composition is continuous because of Lemma \ref{cL_prod}. The first function is continuous because $S\mapsto Se_0$ is continuous and because the mapping $(n,S)\mapsto \zeta_\alpha(n,X(S e_0))$ can be decomposed as

$$
\begin{array}{ccccl}
L^1_{\kappa_0}([x_b,\infty),\mathbb{R}_+)\times\mathbb{R}_+ & \longrightarrow & L^1_{\kappa_0}([x_b,\infty),\mathbb{R}_+)\times W^{g_{\infty}}_\alpha & \longrightarrow & L^1_{-(\hat{\mu}-\mu_0)}(\mathbb{R}_-,\mathbb{R}_+)\\
(n,S) & \longmapsto & (n,X(Se_0)) & \longmapsto & \zeta_\alpha(n,X(S e_0)) \\
\end{array},
$$
so that, by Lemmas \ref{cL_XandTauCont} and \ref{cL_changeVar}, it is also continuous. Therefore $\mathcal{L}^{-1}_{ps}$ is continuous.
\eofproof

\section{Differentiability of $F_1$ and $F_2$}\label{appendix-diffF}

Here we show that the mappings $F_1$ and $F_2$, defining the system of two delay equations, cf. \eqref{sect6-eq1}, \eqref{sect6-eq3} and below, are $C^1$ if certain smoothness and growth conditions on the model ingredients $g, \beta, \mu, \gamma$ and $f$ are satisfied, and the parameter $\mu_0$ characterising the spaces $\X_1$ and $\X_2$ is chosen suitably. Our analysis does not (yet) cover the situation in which growth and reproduction undergo an instantaneous change upon reaching adult size, as described in \cite{DGMNR} (beware and also consult \cite{DGMNR-r}). In \cite{DK}, it is, by way of a related example, shown that in such a situation the solution operators may be $C^1$ even though the mapping defining the delay equation is not. In other words, the Principle of Linearised Stability may very well hold when the sufficient conditions introduced below do not!

We recall the setting:\\

$0<\mu_0<\hat{\mu}$ with $\hat{\mu}$ “defined” by \eqref{weight-ass1},
\begin{equation*}
\begin{array}{ll}

\X_1:= \left\{\phi \in L^1((-\infty,0];\mathbb{R}) : ||\phi||_1 := \displaystyle\int_{-\infty}^0 e^{\mu_0 \sigma}  |\phi(\sigma)| d \sigma < \infty,  \phi \geq 0 \right\}, \\

\X_2:=\left\{\psi \in C((-\infty,0];\mathbb{R}) : ||\psi||_{\infty} := \displaystyle\sup_{\sigma\le 0} \,e^{\mu_0 \sigma} |\psi(\sigma)|  < \infty, \psi \geq 0\right \}.

\end{array}
\end{equation*}
As norm on the Cartesian product $\X_1 \times \X_2$ we choose $||(\phi,\psi)|| = ||\phi||_1 + ||\psi||_{\infty}.$\\
$F_1$ and $F_2\,: \X_1 \times \X_2 \rightarrow \mathbb{R}$ are defined by \\

$F_1 (\phi,\psi) = \langle { H_1(\psi) , \phi } \rangle$ with   $H_1 : \X_2 \rightarrow \X_1’$ with $\X_1’$ the dual space of $\X_1$, represented by a weighted version of $L^{\infty}$, i.e.,
$$
    \X_1’ := \left \{ \theta : [0,\infty)\rightarrow \mathbb{R} : \theta \,\, \text{measurable and}\,\, ||\theta||_{1}’ = \text{ess} \sup \{ |\theta(a)| e^{\mu_0 a} : a \geq 0\} < \infty \right\},
$$
and $H_1(\psi)(a) = \beta(X_{\psi}(0,-a,x_b),\psi(0)) \mathcal{F}_\psi(0,a,x_b);$\\

$F_2(\phi,\psi)(a) = f(\langle{\delta,\psi}\rangle) - \langle {H_2(\psi) , \phi }\rangle$ with $\delta \in H_2'$ defined by $\langle {\delta, \psi} \rangle = \psi(0)$ and 
$H_2 : \X_2 \rightarrow \X_1’,$ $H_2(\psi)(a) = \gamma(X_{\psi}(0,-a,x_b),\psi(0)) \mathcal{F}_\psi(0,a,x_b).$\\

We note that both $F_1$ and $F_2$ are linear in their first argument $\phi$, essentially since we add contributions of individuals. The joyful consequence is that we do not have to worry about the fact that $\X_1$, being a positive cone in an $L^1$ space, has empty interior (side remark: \cite{Ruess} provides an appropriately adapted definition of differentiability). Indeed, when $f$, $H_1$ and $H_2$ are continuously (Fr\'{e}chet) differentiable, then so are $F_1$ and $F_2$ with:
\begin{equation*}
DF_1(\phi_0,\psi_0)(\phi_1,\psi_1) = \langle H_1(\psi_0) , \phi_1 \rangle + \langle DH_1(\psi_0) \psi_1 , \phi_0 \rangle,
\end{equation*}
and
\begin{equation*}
DF_2(\phi_0,\psi_0)(\phi_1,\psi_1) = f’(\langle \delta, \psi_0 \rangle) \langle \delta,\psi_1 \rangle - \langle H_2(\psi_0) , \phi_1 \rangle - \langle DH_2(\psi_0) \psi_1 , \phi_0 \rangle.
\end{equation*}

 Therefore we now concentrate on deriving conditions that guarantee the Fr\'{e}chet differentiability of $H_1$ and $H_2$. In fact we limit our attention to $H_1$, since by copying the assumptions concerning $\beta$ to corresponding assumptions concerning $\gamma$, we cover $H_2$. We first recall that the constructive definition of $H_1$ involves the function $x(\tau)$ with parameters $a$ and $\psi$ defined by the ODE initial value problem 
\begin{equation*}
\begin{aligned}
x'(\tau) = & g(x(\tau),\psi(\tau)), \quad    -a \leq \tau \leq 0,\\
x(-a) = & x_b.
\end{aligned}
\end{equation*}
To simplify slightly the notation without the risk of confusion, we denote here the unique solution by $x(\tau;a,\psi) $ (notice that this was denoted by $X_{\psi}(\tau, -a, x_b)$ in \eqref{X-s-def-2}). Then 
\begin{equation*}
H_1(\psi)(a) := \beta(x(0;a,\psi),\psi(0)) \tilde{f}(0;a,\psi),
\end{equation*}
with 
\begin{equation}\label{directproof}
\tilde{f}(\tau;a,\psi) :=\exp\left\{ - \int_{-a}^{\tau} \mu(x(\eta;a,\psi),\psi(\eta)) \,\ud \eta \right\},\quad -a \leq \tau \leq 0. 
\end{equation}
(Notice that this was denoted by $\mathcal{F}_{\psi}(\tau, -a, x_b)$ in \eqref{F-s-def-2}.)
In order to introduce the candidate for $DH_1(\psi_0) \psi_1$ we need to introduce the variational equation

$$
\begin{aligned}
y'(\tau) = & D_1g(x(\tau;a,\psi_0),\psi_0(\tau)) y(\tau) + D_2g(x(\tau;a,\psi_0),\psi_0(\tau)) \psi_1(\tau), \quad    -a \leq \tau \leq 0,\\
y(-a) = & 0.
\end{aligned}
$$

We shall denote the unique solution of the ODE above by $y(\tau;a,\psi_0,\psi_1).$
Below we shall formulate assumptions that guarantee
$$
(DH_1(\psi_0) \psi_1)(a) = (1) + (2) + (3),
$$
with

$$
\begin{aligned}
(1)	= &  D_2 \beta(x(0;a,\psi_0),\psi_0(0)) \psi_1(0) \tilde{f}(0;a,\psi_0),\\
(2)  = & D_1 \beta(x(0;a,\psi_0),\psi_0(0)) y(0;a,\psi_0,\psi_1) \tilde{f}(0;a,\psi_0), \\
(3)  = &  - \beta(x(0;a,\psi_0),\psi_0(0)) \tilde{f}(0;a,\psi_0) \displaystyle\int_{-a}^0 \left[ D_1 \mu(x(\tau;a,\psi_0),\psi_0(\tau)) y(\tau;a,\psi_0,\psi_1)+ D_2 \mu(x(\tau;a,\psi_0), \psi_0(\tau)) \psi_1(\tau)\right] d \tau.
\end{aligned}
$$

\begin{definition}\label{C1R}
We call $h : \mathbb{R}_+^2 \rightarrow \mathbb{R}_+$ “regular enough” with parameters $L$ and $C$ if

\renewcommand{\labelenumi}{(\roman{enumi})}
\begin{enumerate}
	\item	$h$ is a $C^1$ map .
	\item	$h$ is globally Lipschitz continuous with constant $L$, i.e., 
	$$
	|h(x+\xi, S + \sigma) - h(x,S)| \leq L ( |\xi| + |\sigma| ).
	$$
	Note that, as a consequence, both partial derivatives are uniformly bounded by $L$.
	\item  With R defined by
	$$
	R(x,S,\xi,\sigma)= h(x+\xi,S+\sigma) - h(x,S) - D_1h(x,S) \xi - D_2h(x,S) \sigma
	$$
	there exists a constant $C>0$ such that, uniformly for $(x,S)$ in $\mathbb{R}_{+}^2,$
	$$
	     | R(x,S,\xi,\sigma)| \leq C ( |\xi| + |\sigma| ) ^2.
	$$
\end{enumerate}
\end{definition}
Notice that $h$ is regular enough if hypotheses H1$_h$, H4$_h$ and H5$_h$ hold. We assume that $g$ is regular enough with parameters $L_1$ and $C_1$. 

In the following $\psi_0$ denotes an element of $\X_2$. $\psi_1$, on the other hand, denotes an element of norm one in the corresponding weighted space of continuous functions that may take negative values, but $\psi_1$ should be such that the sum $\psi_0 + \varepsilon \psi_1$ belongs to $\X_2$ for small positive values of $\varepsilon$. In order to simplify, we slightly abuse notation by writing
$$
x(\tau;a,\varepsilon) \quad \text{to denote} \quad x(\tau;a,\psi_0 + \varepsilon \psi_1)
$$
and, similarly,
$$
    \tilde{f}(\tau;a,\varepsilon) \quad \text{to denote} \quad \tilde{f}(\tau;a,\psi_0 + \varepsilon \psi_1).
$$
\begin{lemma}\label{C1}
For $a \geq 0$ and $-a \leq \tau \leq 0$, the inequality
\begin{equation*}
    |x(\tau;a,\varepsilon) - x(\tau;a,0)| \leq \varepsilon L_1 (a+\tau) e^{L_1(a+\tau)} 
\sup\{|\psi_1(\sigma)| : -a \leq \sigma \leq \tau\}
\end{equation*}
holds.
\end{lemma}
{\bf Proof.} 
$$
x(\tau;a,\varepsilon) = x_b + \int_{-a}^{\tau} g(x(\sigma;a,\varepsilon), \psi_0(\sigma) + \varepsilon \psi_1(\sigma)) d \sigma
$$
Hence
$$
x(\tau;a,\varepsilon) - x(\tau;a,0) = \int_{-a}^{\tau} \bigg(g(x(\sigma;a,\varepsilon), \psi_0(\sigma) + \varepsilon \psi_1(\sigma))
	-     g(x(\sigma;a,0), \psi_0(\sigma)) \bigg) d \sigma.
$$
Now use property (ii) of Definition \ref{C1R} to obtain
$$
|x(\tau;a,\varepsilon) - x(\tau;a,0)| \leq L_1 \int_{-a}^{\tau} |x(\sigma;a,\epsilon) - x(\sigma;a,0)| d \sigma + \varepsilon L_1 \int_{-a}^{\tau} |\psi_1(\sigma)| d \sigma.
$$
Next apply Gr\"{o}nwall's inequality. 

\eofproof

By applying Gr\"{o}nwall's inequality directly to the equation for $y$ we obtain the following lemma.
\begin{lemma}\label{C2}
	For $a \geq 0$ and $-a \leq \tau \leq 0$ the inequality
	$$
	|y(\tau;a,\psi_0,\psi_1)| \leq  L_1 (a+\tau) e^{L_1(a+\tau)} \sup\{|\psi_1(\sigma)| : -a \leq \sigma \leq \tau\}
	$$
	holds.
\end{lemma}
Our aim is to derive estimates for $z(\tau;a,\varepsilon)$ defined by
$$
z(\tau;a,\varepsilon) = \frac{1}{\varepsilon} \big(x(\tau;a,\varepsilon) - x(\tau;a,0)\big) - y(\tau;a,\psi_0,\psi_1).
$$
By combining Lemmas \ref{C1} and \ref{C2} we obtain the following estimate 
$$
   |z(\tau;a,\varepsilon)| \leq  2 L_1 (a+\tau) e^{L_1(a+\tau)} \sup\{|\psi_1(\sigma)| : -a \leq \sigma \leq \tau \},
$$
but we can bootstrap by first observing that $z$ satisfies (with $R$ as introduced in Definition \ref{C1R} but now for $h=g$)

$$
\begin{aligned}
z’(\tau) = & D_1g(x(\tau;a,\psi_0),\psi_0(\tau)) z(\tau) + 
\frac{1}{\varepsilon} R(x(\tau;a,0),\psi_0(\tau),x(\tau;a,\varepsilon)-x(\tau;a,0),\varepsilon \psi_1(\tau))\\
z(-a) = & 0
\end{aligned}
$$
 \begin{lemma}\label{C3}
 	For $a \geq 0$ and $-a \leq \tau \leq 0$ the inequality
 	$$
    |z(\tau;a,\varepsilon)| \leq \varepsilon C_1 (a+\tau) e^{L_1(a+\tau)} \big(L_1(a+\tau) e^{L_1(a+\tau)} + 1 \big)^2 
 \sup\{|\psi_1(\sigma)|^2 : -a \leq \sigma \leq \tau\}
 	$$
 	holds.
 \end{lemma}
{\bf Proof.}
The combination of the differential equation for $z$ and the initial condition imply

$$
z(\tau;a,\varepsilon) = \frac{1}{\varepsilon}\int_{-a}^{\tau} \exp\left\{\int_{\sigma}^{\tau} D_1g(x(\eta;a,0),\psi_0(\eta)) \,\ud \eta\right\}
R(x(\sigma;a,0),\psi_0(\sigma),x(\sigma;a,\varepsilon) - x(\sigma;a,0),\varepsilon \psi_1(\sigma))\,\ud \sigma.
$$

Using Assumption (iii) of Definition \ref{C1R}, Lemma \ref{C1} and the fact that $D_1g$ is bounded by $L_1$ we obtain the estimate stated in the lemma.
\eofproof
\begin{corollary}\label{C4}
	For given $a > 0$, $y( \cdot ;a,\psi_0,\psi_1)$ is the Fr\'{e}chet derivative of $x( \cdot ;a,\psi)$ with respect to $\psi$, taken in $\psi_0$ and acting on $\psi_1$, when these functions, including $\psi$, are considered as elements of $C([-a,0];\mathbb{R})$ equipped with the supremum norm.
\end{corollary}
   Our interest, however, is in $x(0;a,\psi)$ as a function of $a$, with $a$ ranging in $[0,\infty)$. When 
$||\psi_1|| = 1$,  $\sup\{|\psi_1(\sigma)|^2 : -a \leq \sigma \leq \tau\}$ can grow like $e^{ 2 \mu_0 a }$ and we shall need to cope with that growth. The more worrisome feature is the growth of the multiplicative factor $a e^{L_1a}$ for large $a$. The factor takes this form because we replaced $D_1 g$ by $L_1$ to obtain estimates. For large $a$, size may be large as well and $D_1 g$ might then be much smaller than $L_1$. There are many ways in which we can formalize this idea. We have chosen a somewhat drastic version.\\
\\
\textbf{Assumption} (cf. \eqref{growth-assum} i.e. H2$_g$ and H3$_g$; and H$_{g_{\infty}}$) \\
\renewcommand{\labelenumi}{(\roman{enumi})}
\begin{enumerate}
	\item 	There exists $g_{\max}\ge g_{\min} > 0$ such that $g_{\max}\ge g(x,S) \geq g_{\min}$ for $(x,S) \in [x_b, \infty) \times \mathbb{R}_{+}$.
	\item	There exists $\bar{x} > x_b$ and $g_{\infty} \geq g_{\min}$ such that $g(x,S) = g_{\infty}$ for $x \geq \bar{x}.$
\end{enumerate}
As an immediate consequence we have
\begin{lemma}\label{C6}
	Let $\bar{a}$ be the maximum time it can take an individual to grow from birth size $x_b$ to the size $\bar{x}$ at which the growth rate becomes the constant $g_{\infty}$, i.e.,
	$$
	\bar{a} := \frac{\bar{x} - x_b}{g_{min}}.
	$$
	Then for $a > \bar{a}$ and $\tau \geq \bar{a} - a$,
    \renewcommand{\labelenumi}{(\roman{enumi})}
    \begin{enumerate}
		\item $x(\tau;a,\varepsilon) = x(\bar{a} - a;a,\varepsilon) + g_{\infty} (\tau - \bar{a} + a),$\\
		\item	$y(\tau;a,\psi_0,\psi_1) = y(\bar{a} - a;a,\psi_0,\psi_1)$,\\
		\item 	$z(\tau;a,\varepsilon) = z(\bar{a} - a;a,\varepsilon)$.
	\end{enumerate}	
\end{lemma}
\begin{corollary}\label{C7}
	For all $a \geq 0,$
	   $$
	   |z(0;a,\varepsilon)| \leq \varepsilon \, C_1 \bar{a} e^{L_1 \bar{a}} \left( L_1 \bar{a} e^{L_1 \bar{a}} + 1 \right)^2 
	\sup\{|\psi_1(\sigma)|^2 : -a \leq \sigma \leq \min (-a+\bar{a}, 0 ) \}.
	$$
\end{corollary}
We now assume that both $\beta$ and $\mu$ are “regular enough” too, with parameters, respectively, $(L_2,C_2)$  and $(L_3,C_3)$ (notice that $\beta$ and $\mu$ are “regular enough” if hypothesis H1$_{\beta}$, H4$_{\beta},$ H5$_{\beta,} $H1$_{\mu}$, H4$_{\mu}$ and H5$_{\mu}$ hold).
\begin{lemma}\label{C8}
$$
   \left| \beta(x(0;a,\varepsilon), \psi_0(0) + \varepsilon \psi_1(0)) - \beta(x(0;a,0),\psi_0(0)) - 
\varepsilon D_1\beta(x(0;a,0),\psi_0(0)) y(0;a,\psi_0,\psi_1) - \varepsilon D_2 \beta(x(0;a,0),\psi_0(0)) \psi_1(0)\right |
$$
$$
\leq C_4\, \varepsilon^2 \, \sup\{|\psi_1(\sigma)|^2 : -a \leq \sigma \leq 0 \}.
$$	
\end{lemma}
{\bf Proof.} First, observe that
$$
   x(0;a,\varepsilon) = x(0;a,0) + \varepsilon y(0;a,\psi_0,\psi_1) + \varepsilon z(0;a,\varepsilon).
$$
With $R$ as introduced in Definition 1, but now for $h= \beta$, the left hand side of the inequality can be written as
$$
   \left | R\Big(x(0;a,0), \psi_0(0), x(0;a,\varepsilon) - x(0;a,0) ,  \varepsilon \psi_1(0)\Big) +D_1 \beta(x(0;a,0),\psi_0(0)) \varepsilon z(0;a,\varepsilon)\right |.
$$
Now note that
$$
\left | R\left (x(0;a,0), \psi_0(0), x(0;a,\varepsilon) - x(0;a,0), \varepsilon \psi_1(0) \right ) \right | \leq C_2 \left( \left|x(0;a,\varepsilon) - x(0;a,0) \right | + \varepsilon |\psi_1(0)| \right )^2,
$$
and hence, by Lemma \ref{C1}
$$
\left| R(x(0;a,0), \psi_0(0), x(0;a,\varepsilon) - x(0;a,0), \varepsilon \psi_1(0)) \right | \leq C_2 \varepsilon^2 \left(L_1 \bar{a} e^{L_1 \bar{a}} + 1 \right )^2 \sup\{|\psi_1(\sigma)|^2 : -a \leq \sigma \leq 0 \}.
$$
Since $|D_1 \beta| \leq L_2$ we obtain from Corollary \ref{C7} the estimate 
$$
|D_1 \beta(x(0;a,0), \psi_0(0)) \varepsilon z(0;a,\varepsilon)| \leq L_2 \varepsilon^2 C_1 \bar{a} e^{L_1 \bar{a}} \left ( L_1 \bar{a} e^{L_1 \bar{a}} + 1 \right )^2 \sup\{|\psi_1(\sigma)|^2 : -a \leq \sigma \leq \min( -a+\bar{a},0)\}.
$$
Combination of these two inequalities leads, for a suitable choice of $C_4$, to the statement of the lemma.
\eofproof\\
In exactly the same way one proves
\begin{lemma}\label{C9}
	$$
\left | \tilde{f}(0;a,\varepsilon) - \tilde{f}(0;a,0) - \varepsilon \tilde{f}(0;a,0) \int_{-a}^0 D_1 \mu(x(\tau;a,\psi_0),\psi_0(\tau)) y(\tau;a,\psi_0,\psi_1)
	+ D_2 \mu(x(\tau;a,\psi_0),\psi_0(\tau)) \psi_1(\tau)\,\, d \tau \right | 
	$$
	$$
	\leq C_5 \varepsilon^2 \sup\{|\psi_1(\sigma)|^2 : -a \leq \sigma \leq \min(-a+\bar{a},0)\}
	$$	
\end{lemma}
\begin{theorem}
	Assume that $g$, $\beta$ and $\mu$ are ‘regular enough’, that \eqref{weight-ass1} and Assumption 5 (equivalently \eqref{growth-assum} and H$_{g_{\infty}}$) hold, that $\beta$, $\mu$, as well as their first order partial derivatives, are bounded and that $3 \mu_0 < \hat{\mu}$. \\
	Then $H_1 : \X_2 \rightarrow \X_1’$ is continuously Fr\'{e}chet differentiable with derivative given by  
	$$
	\begin{array}{ll}
	(DH_1(\psi_0)\psi_1)(a) \\
	= \tilde{f}(0;a,\psi_0) \bigg( D_2 \beta(x(0;a,\psi_0),\psi_0(0)) \psi_1(0) + D_1 \beta(x(0;a,\psi_0),\psi_0(0)) y(0;a,\psi_0,\psi_1)  -\\ \beta(x(0;a,\psi_0),\psi_0(0)) \int_{-a}^0 D_1 \mu(x(\tau;a,\psi_0),\psi_0(\tau)) y(\tau;a,\psi_0,\psi_1) + D_2 \mu(x(\tau;a,\psi_0),\psi_0(\tau)) \psi_1(\tau) \,\,  d \tau  \bigg) 
	\end{array}
	$$
\end{theorem}
{\bf Proof.} We first show that $DH_1(\psi_0)$ is a bounded linear operator from $\X_2$ to $\X_1’$, depending continuously on $\psi_0 \in \X_2$. By combining Lemma \ref{C6}.ii and Lemma \ref{C2} we obtain

$$
|y(\tau;a,\psi_0,\psi_1)| \leq L_1 \bar{a} e^{L_1 \bar{a}} \sup\{|\psi_1(\sigma)| : -a \leq \sigma \leq \min(-a+\bar{a},0) \}
$$

Since $|\psi_1(\sigma)| e^{\mu_0 \sigma} \leq ||\psi_1||_{\infty}$ and $e^{- \mu_0 \sigma} \leq e^{\mu_0 a}$ for $\sigma \geq -a$, it follows that for any $\tau \leq 0$

$$
|y(\tau;a,\psi_0,\psi_1)| \leq C_6 e^{\mu_0 a} ||\psi_1||_{\infty}.
$$

More directly it follows that $|\psi_1(\tau)| \leq e^{\mu_0 a} ||\psi_1||_{\infty}$ for $-a \leq \tau \leq 0$.
If we multiply the expression for $(D H_1(\psi_0)\psi_1)(a)$ by $e^{\mu_0 a},$ use these estimates, and take the supremum with respect to $a$, we obtain $||\psi_1||_{\infty}$ multiplied by a scalar factor that is finite since $ \sup \{ a e^{2 \mu_0 a} \tilde{f}(0;a,\psi_0) \} $ is finite by \eqref{weight-ass1} and the condition that $3 \mu_0 < \hat{\mu}.$

We now proceed to prove the continuity of the differential. From the formulas for $D_1 F$ and $D_2 F$ in terms of $H_1, H_2$ and $f,$ it is clear that we only need to check the continuity of
\begin{equation*}
\begin{aligned}
&\X_2 & \longrightarrow &\quad \mathcal{BL}(\X_2, \X_1'),\\
&\psi_0 & \longrightarrow &\quad DH_i (\psi_0),
\end{aligned}
\end{equation*}
for $i = 1,2$. (Above $\mathcal{BL}(\X_2, \X_1')$ stands for the Banach space of bounded linear operators from $\X_2$ to $\X_1'$.)

The notation introduced above Lemma \ref{C1} is not very convenient here. So we go back to denote by $x(\tau;a,\psi)$ the size at time $\tau$ of an individual of age $a$ which has experienced a resource level $\psi$. Recall that the same was called $X_{\psi}(\tau, -a, x_b)$ in Section \ref{sectPDEformulation} (see \eqref{X-s-def-2}); also $\tilde{f} (\tau;a,\psi)$ is what we earlier denoted by $\mathcal{F}_{\psi}(\tau,-a,x_b)$ (see \eqref{F-s-def-2}). With this, Lemma \ref{C1} can be reformulated as 
\begin{lemma}\label{reform}
For $a \geq 0$ and $-a \leq \tau \leq 0$ the inequality
$$
|x(\tau;a,\psi_0 + \tilde{\psi_0}) - x(\tau;a,\psi_0)| \leq  L_1 (a+\tau) e^{L_1(a+\tau)} 
\sup_{-a \leq \sigma \leq 0}|\tilde{\psi_0}(\sigma)| 
$$
holds.
	\end{lemma}
Recall the three term decomposition of $DH_1$ introduced above Definition C.1. Let us first focus on the term (1). We need to show that (recall that only the value of $\psi_1$ at $0$ has an influence on the first term)
\begin{equation}\label{cont_first}
\begin{array}{cl}
\displaystyle\sup_{|\psi_1(0)|\leq 1} \, \sup_{a \geq 0} \Big \{ e^{\mu_0 a} |D_2 \beta \big(x(0;a,\psi_0 + \tilde{\psi}_0),\psi_0(0)+\tilde{\psi_0}(0)\big) \psi_1(0)  \tilde{f}(0;a,\psi_0 + \tilde{\psi}_0) \\- D_2 \beta \big(x(0;a,\psi_0) ,\psi_0 (0) \big)\psi_1(0) \tilde{f}(0;a,\psi_0)| \Big \} \longrightarrow 0
\end{array}
\end{equation}
when $|| \tilde{\psi}_0 ||_{\X_2} \rightarrow 0.$
By replacing $\psi_1(0)$ by $1$ we take care of the first sup. Our next step is to show that we can restrict to a bounded set for the variable $a,$ since the whole expression converges to zero for $a \rightarrow \infty$, uniformly in $\tilde{\psi}_0.$ We assume\\
\textbf{A1} \,\,\, $D_2 \beta$ is uniformly bounded.\\
Recalling (\ref{weight-ass1}), we replace in (\ref{cont_first})  $e^{\mu_0 a}$ by $e^{(\mu_0 - \hat{\mu}) a}$ and multiply both $\tilde{f}(0;a,\psi_0 + \tilde{\psi}_0)$ and $\tilde{f}(0;a,\psi_0 )$ by $e^{\hat{\mu} a}.$ The factor between vertical bars is bounded, uniformly in $\tilde{\psi}_0;$ whereas the factor $e^{(\mu_0 - \hat{\mu}) a}$ tends to $0$ as $a$ goes to $\infty.$ Hence, for any $\varepsilon>0$ there exists $\tilde{a}(\varepsilon)$ such that the product can be bounded by $\varepsilon$ for $0 \leq a \leq \tilde{a}(\varepsilon)$ by making $||\tilde{\psi}_0||_{\X_2}$ sufficiently small. We also assume\\
\textbf{A2} \,\,\, $D_2 \beta$ and $\mu$ are globally Lipschitz continuous.\\
We show  that $\psi_0 \mapsto D_2 \beta\big(x(0;a, \psi_0),\psi_0(0)\big)$ and $\psi_0 \mapsto \tilde{f}(0;a, \psi_0)$ are continuous as maps from $\X_2$ to $\mathbb{R},$ uniformly for $a$ in compact sets. The result next follows from the standard result that the product of two continuous functions is continuous. Note that the function $e^{\mu_0 a}$ is bounded on $[0, \tilde{a}(\varepsilon)].$\\
The continuity of $\psi_0 \mapsto D_2 \beta \big( x(0;a,\psi_0), \psi_0(0) \big)$ is a consequence of the Lipschitz continuity of $D_2 \beta$ and Lemma \ref{reform}.\\
Concerning  $\psi_0 \mapsto \tilde{f}(0;a,\psi_0),$ a stronger result is already available, see Lemma \ref{C9}, under aditional assumptions on $\mu.$  A direct proof starts from
\begin{equation}\label{ftilde}
\tilde{f}(\tau;a,\psi):= \exp\left\{- \int_{-a}^{\tau}\mu\big(x(\eta;a,\psi), \psi(\eta) \big) \,\ud \eta\right\}
\end{equation}
(see \ref{directproof}). We write
\begin{align*}
\tilde{f}(\tau;a,\psi_0 + \tilde{\psi}_0)& -\tilde{f}(\tau;a,\psi_0) = \\ 
& \tilde{f}(\tau;a,\psi_0) \left\{\exp  \left\{\int_{-a}^{\tau} \left(\mu\left(x(\eta;a,\psi_0),\psi_0(\eta) \right) - \mu\left( x(\eta;a,\psi_0 + \tilde{\psi}_0), \psi_0(\eta) + \tilde{\psi}_0(\eta) \right)  \right)\ud \eta \right\}  - 1\right\}
\end{align*}
and use
\begin{enumerate}
	\item $|e^y - 1| < 2 |y|$ for $|y|$ small.
	\item Lipschitz continuity of $\mu.$
	\item Lemma \ref{reform}.
\end{enumerate}
This concludes the analysis of the term (1).
We now focus on the term (2). In the assumptions A1 and A2 above we change $D_2$ to $D_1$. The factor $\psi_1(0)$ that we had before is now replaced by $y(0;a,\psi_0, \psi_1)$ to which the estimate of Lemma \ref{C2} and, in addition, the observation of Lemma \ref{C6} (ii) apply. This leads to 
$$
|y(0;a, \psi_0, \psi_1) | \leq L_1 \bar{a} e^{L_1 \bar{a}} \sup \{|\psi_1(\sigma)|: -a \leq \sigma \leq \min\{\bar{a} -a,0 \} \}
$$
We note that \, $||\psi_1 || < 1$ if and only if $|\psi_1 (\sigma)| < e^{\mu_0 a}, \sigma \leq 0.$ So the factor $\sup\{|\psi_1(\sigma)|:-a \leq \sigma \leq \min \{ \bar{a}-a,0\}\}$ can be estimated by $e^{\mu_0 a}.$ Compared to the term (1), this yields an extra factor $e^{\mu_0 a},$ but as long as $2 \mu_0 - \hat{\mu} < 0,$ all the preceding arguments work just as well. This establishes the appropriate continuity property of the term (2).\\
In order to deal with the term (3), we assume\\
\textbf{A3} \,\,\, $D_1 \mu$ and $D_2 \mu$ are globally Lipschitz continuous.\\
The estimates
\begin{align*}
|y(\tau;a,\psi_0,\psi_1)| & \leq L_1 \bar{a} e^{L_1 \bar{a}} e^{\mu_0 a},\\
|\psi_1(\tau)| & \leq e^{\mu_0 a},\quad   -a \leq \tau \leq 0
\end{align*}
follow from $||\psi_1|| \leq 1 $ exactly as sketched above.\\
Thus the continuity of the term (3) can be shown in the same manner as the continuity of the term (2).

   Finally, it remains to prove that $DH_1(\psi_0)$ is indeed the derivative in $\psi_0.$
In our simplified notation we have 
$$
H_1(\psi_0 + \varepsilon \psi_1)(a) = \beta (x(0;a,\varepsilon),\psi_0(0) + \varepsilon \psi_1(0)) \tilde{f}(0;a,\varepsilon)
$$
and hence
\\
$$
\begin{aligned}
& H_1(\psi_0 + \varepsilon \psi_1)(a) - H_1(\psi_0)(a) \\
= & \Big(\beta(x(0;a,\varepsilon), \psi_0(0) + \varepsilon \psi_1(0)) - \beta(x(0;a,0),\psi_0(0))\Big) \tilde{f}(0;a,\varepsilon) + \beta(x(0;a,0),\psi_0(0)) \big(\tilde{f}(0;a,\varepsilon) - \tilde{f}(0;a,0)\big)  \\
= & \bigg(\beta(x(0;a,\varepsilon), \psi_0(0) + \varepsilon \psi_1(0)) - \beta(x(0;a,0),\psi_0(0))\\ 
& -\varepsilon D_1\beta(x(0;a,0),\psi_0(0)) y(0;a,\psi_0,\psi_1) - \varepsilon D_2\beta(x(0;a,0),\psi_0(0)) \psi_1(0)\bigg) \tilde{f}(0;a,\varepsilon)\\  
& + \beta(x(0;a,0),\psi_0(0)) \bigg( \tilde{f}(0;a,\varepsilon)   - \tilde{f}(0;a,0) \\
& - \varepsilon \tilde{f}(0;a,0) \int_{-a}^0 D_1 \mu(x(\tau;a,\psi_0),\psi_0(\tau)) y(\tau;a,\psi_0,\psi_1) + D_2 \mu(x(\tau;a,\psi_0),\psi_0(\tau)) \psi_1(\tau) \, d \tau \bigg)\\ 
& +	\bigg(\varepsilon D_1\beta(x(0;a,0),\psi_0(0)) y(0;a,\psi_0,\psi_1) + \varepsilon D_2 \beta(x(0;a,0),\psi_0(0)) \psi_1(0)\bigg) \tilde{f}(0;a,\varepsilon)\\
& -	\beta(x(0;a,0),\psi_0(0)) \varepsilon \tilde{f}(0;a,0) \int_{-a}^0 D_1 \mu(x(\tau;a,\psi_0),\psi_0(\tau)) y(\tau;a,\psi_0,\psi_1)+	D_2\mu(x(\tau;a,\psi_0),\psi_0(\tau)) \psi_1(\tau)\,\, d\tau .
\end{aligned}
$$
\\
\\
If we now divide this identity by $\varepsilon$ and let $\varepsilon$ go to zero, the last two terms converge, pointwise in $a$, to $(D H_1(\psi_0)\psi_1)(a),$ while the first two terms converge to zero on account of, respectively, Lemma \ref{C8} and Lemma \ref{C9}. If we multiply by $e^{\mu_0 a}$ and take the supremum with respect to $a$, the convergence still holds, since uniformly for $\psi_1$ of norm $1$ the estimate
$$
\sup\left\{|\psi_1(\sigma)|^2 : -a \leq \sigma \leq \min(-a+\bar{a},0)\right \} \leq e^{2 \mu_0 a}
$$
holds and the $\tilde{f}$ decays, because of \eqref{weight-ass1}, sufficiently fast to let the integral converge.
\eofproof\\
\\
We are now ready to state the main result.

\begin{theorem}\label{theoremC12}
Assume that H4$_f$,  H3$_g$ , H$_s$, H$_{g_\infty}$ as well as H1$_h$, H2$_h$, H4$_h$, H5$_h$ for $h=g, \beta, \gamma,\mu$ hold.
Then $F_1$ and $F_2$ are continuously Fr\'{e}chet differentiable maps from $\mathcal{X}_1\times\mathcal{X}_2$ to $\mathbb{R}$,  when the weight $\mu_0$ satisfies $3\mu_0< \hat{\mu}$.
\end{theorem}

\section*{Acknowledgments}
We thank the International Centre for Mathematical Sciences for financial 
support we received from the Research in Groups program during our stay at Edinburgh in July 2017; and the Spanish research projects MTM2014-52402-C3-2P and MTM2017-84214-C2-2P.

\end{document}